\documentclass[12pt,final]{amsart}
\usepackage[final]{graphicx}
\usepackage{amssymb}
\usepackage{alltt}
\usepackage{verbatim} 

\newcommand{\regfigure}[2]{\includegraphics[height=#2in,
keepaspectratio=true]{#1.eps}}
\newcommand{\widefigure}[3]{\includegraphics[height=#2in,
width=#3in]{#1.eps}}

\newtheorem{thm}{Theorem}
\newtheorem{lem}[thm]{Lemma}
\newtheorem{cor}[thm]{Corollary}
\newtheorem*{thm*}{Theorem}
\newtheorem*{lem*}{Lemma}
\newtheorem*{cor*}{Corollary}

\newtheorem*{prop}{Proposition}
\theoremstyle{definition}
\newtheorem*{defn}{Definition}
\newtheorem*{defns}{Definitions}
\newtheorem{example}{Example}
\newtheorem{exer}{Exercise}
\newtheorem*{rem}{Remark}

\newcommand{\CC}{\mathbb{C}}

\newcommand{\RR}{\mathbb{R}}
\newcommand{\ZZ}{\mathbb{Z}}

\newcommand{\ba}{\mathbf{a}}
\newcommand{\balpha}{\mathbf{\alpha}}

\newcommand{\bC}{\mathbf{C}}
\newcommand{\bGamma}{\mathbf{\Gamma}}
\newcommand{\bD}{\mathbf{D}}
\newcommand{\be}{\mathbf{e}}
\newcommand{\bE}{\mathbf{E}}
\newcommand{\bbf}{\mathbf{f}}
\newcommand{\bg}{\mathbf{g}}
\newcommand{\bH}{\mathbf{H}}
\newcommand{\bI}{\mathbf{I}}
\newcommand{\bl}{\mathbf{l}}
\newcommand{\bL}{\mathbf{L}}
\newcommand{\bLambda}{\mathbf{\Lambda}}

\newcommand{\bp}{\mathbf{p}}
\newcommand{\bP}{\mathbf{P}}
\newcommand{\bPi}{\mathbf{\Pi}}
\newcommand{\bq}{\mathbf{q}}
\newcommand{\br}{\mathbf{r}}
\newcommand{\bu}{\mathbf{u}}
\newcommand{\bU}{\mathbf{U}}
\newcommand{\bV}{\mathbf{V}}
\newcommand{\bv}{\mathbf{v}}
\newcommand{\bw}{\mathbf{w}}
\newcommand{\bx}{\mathbf{x}}
\newcommand{\by}{\mathbf{y}}
\newcommand{\bz}{\mathbf{z}}

\newcommand{\Bcal}{\mathcal{B}}
\newcommand{\Co}{\mathcal{C}}
\newcommand{\Ecal}{\mathcal{E}}
\newcommand{\bEcal}{\mathbf{\mathcal{E}}}
\newcommand{\Lop}{\mathcal{L}}

\newcommand{\Po}{\mathcal{P}}
\newcommand{\Pcal}{\mathcal{P}}
\newcommand{\bPo}{\mathbf{\mathcal{P}}}
\newcommand{\Vcal}{\mathcal{V}}
\newcommand{\Xcal}{\mathcal{X}}
\newcommand{\Ycal}{\mathcal{Y}}

\newcommand{\hD}{\hat D}
\newcommand{\hbD}{\mathbf{\hat D}}
\newcommand{\hf}{\hat f}
\newcommand{\hg}{\hat g}
\newcommand{\hbLambda}{\mathbf{\hat \Lambda}}
\newcommand{\hM}{\hat M}

\newcommand{\hT}{\hat T}
\newcommand{\hu}{\hat u}
\newcommand{\hbu}{\mathbf{\hat u}}
\newcommand{\hbU}{\mathbf{\hat U}}
\newcommand{\hU}{\hat U}
\newcommand{\hV}{\hat V}
\newcommand{\hbV}{\mathbf{\hat V}}

\newcommand{\tT}{\tilde T}
\newcommand{\uf}{\underline{f}}

\newcommand{\trefoneline}[2]{\ensuremath{\left[\begin{array}{c|c|c} #1 & \dots & #2  \end{array}\right]}}
\newcommand{\trefonelinefive}[5]{\ensuremath{\left[\begin{array}{c|c|c|c|c} #1 & #2 & #3 & #4 & #5 \end{array}\right]}}

\newcommand{\cond}{\operatorname{cond}}
\newcommand{\diag}{\operatorname{diag}}
\newcommand{\ddt}[1]{\ensuremath{\frac{d#1}{dt}}}
\newcommand{\eps}{\epsilon}

\renewcommand{\Im}{\operatorname{Im}}
\newcommand{\ip}[2]{\ensuremath{\left<#1,#2\right>}}
\newcommand{\lam}{\lambda}
\newcommand{\Matlab}{\textsc{Matlab}}
\newcommand{\mtt}{\texttt}
\renewcommand{\Re}{\operatorname{Re}}

\begin{document}

\title[Chebyshev collocation on ODEs and DDEs]{Chebyshev collocation for linear, periodic \\ ordinary and delay differential equations:\\ \emph{a posteriori} estimates}

\author{Ed Bueler}
\address{Dept. of Mathematical Sciences\\University of Alaska, Fairbanks}
\email{ffelb@uaf.edu}

\thanks{\today.  Supported in part by NSF Grant No. 0114500.}

\begin{abstract}  We present a Chebyshev collocation method for linear ODE and DDE problems.  Theorem \ref{thm:apostivp} in section \ref{sect:apost} gives an \emph{a posteriori} estimate for the accuracy of the approximate solution of a scalar ODE initial value problem.  Examples of the success of the estimate are given.  For linear, periodic DDEs with integer delays we define and discuss the monodromy operator $U$ in section \ref{sect:DDEop}.  Our main goal is reliable estimation of the stability of such DDEs.  In section \ref{sect:DDEapprox} we prove theorem \ref{thm:main} which gives \emph{a posteriori} estimates for eigenvalues of $U$, our main result.  Theorem \ref{thm:main} is based on theorem \ref{thm:BFH}, a generalization to operators on Hilbert spaces of the Bauer-Fike theorem \cite{BauerFike} for (matrix) eigenvalue perturbation problems.  Section \ref{sect:systems} describes the generalization of these results to systems of DDEs and an example is given in section \ref{sect:Mathieu}.  The computation of good bounds on ODE fundamental solutions is an important technical issue and an \emph{a posteriori} method for such bounds is given.  An additional technical issue is addressed in section \ref{sect:polytechs}, namely the evaluation of polynomials and of the $L^\infty$ norms of analytic functions.  Generalization to the non-integer delays case is considered in section \ref{sect:taunoteqperiod}.
\end{abstract}

\setcounter{tocdepth}{1}

\maketitle
\scriptsize\tableofcontents\normalsize

\section{Introduction} \label{sect:intro}

\begin{figure}[ht]
\regfigure{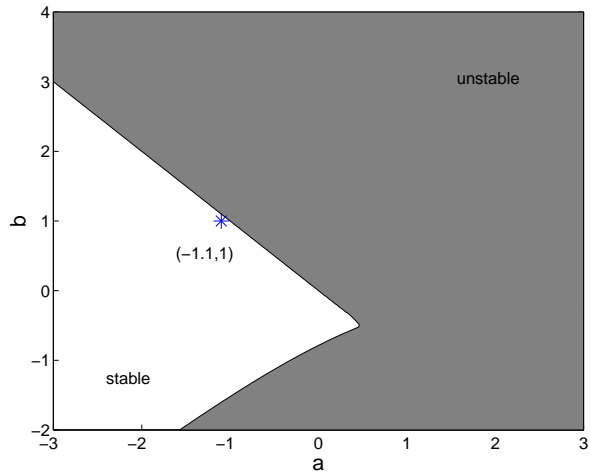}{2.5}
\caption{A stability chart for equation \eqref{introex} with a (numerically-determined) stable point $(a,b)=(-1.1,1)$ indicated.} \label{fig:introchart}
\end{figure}

Consider the following delay differential equation (DDE):
\begin{equation}\label{introex}
    \dot y(t) = a y(t) + \left(b + \sin(3 \pi t)\right) y(t-2).
\end{equation}
    \begin{quote}\textbf{Question:}  For what values of $a,b\in \RR$ is equation \eqref{introex} stable in the sense that all solutions $y(t)$ go to zero as $t\to \infty$?\end{quote}
For $a,b\in[-3,3]\times[-2,4]$, the answer is given by figure \ref{fig:introchart}, a \emph{stability chart}.

Pictures like figure \ref{fig:introchart} of the stability of linear DDE are useful in applications.  Examples include engineering design and control problems for systems described by DDEs \cite{BMBAS,Stepan}, but also applications in biology \cite{MacDonald}, for instance.  

For constant-coefficient cases the stability can be determined ``by hand'' using complex variable techniques but, in general, the stability of linear DDE systems must be determined approximately.  If the Floquet transition matrix for a linear, periodic ordinary differential equation (ODE) associated to a linear, periodic DDE with integer delays can be found exactly, then complex variable techniques will determine stability of the DDE (see section 8.3 of \cite{HaleLunel}).  In general, however, ODE Floquet matrices must themselves be approximated.  The techniques of this paper directly approximate DDE stability by directly approximating the DDE \emph{monodromy operator} (section \ref{sect:DDEop}) associated to equation \eqref{introex}.

Figure \ref{fig:introchart} was produced pixel-by-pixel by approximating the monodromy operator for the parameter values $a,b$ associated to each pixel in the chart, and (numerically) computing the largest eigenvalue of the approximating matrix, a \emph{monodromy matrix}.  The approximation method can be described as Chebyshev (spectral) collocation.

The first question about this procedure, mathematically, is how accurately one can solve initial value problems for DDE \eqref{introex}?  That is, how accurately can one perform the ``method-of-steps'' \cite{HaleLunel} for the DDE?  This paper proves \emph{a posteriori} estimates for the accuracy of the collocation algorithm on such problems.

The next question, mathematically, is about the quality of the eigenvalues of the monodromy matrix as estimates of the eigenvalues of the monodromy operator.  The monodromy operator is compact and thus its spectrum is entirely eigenvalue spectrum (appendix B).  In this paper we combine the results of section \ref{sect:apost} with an eigenvalue perturbation theorem for operators on Hilbert spaces (theorem \ref{thm:BFH} in section \ref{sect:DDEapprox}) to give \emph{a posteriori} estimates for the distance between a nonzero eigenvalue of the monodromy operator and the nearest eigenvalue of the monodromy matrix.

\begin{figure}[ht]
\regfigure{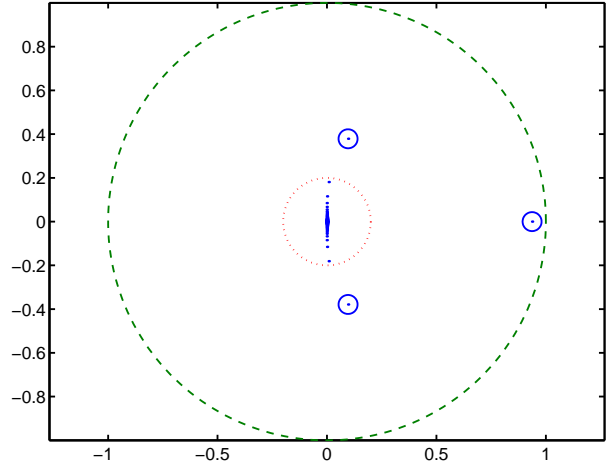}{2.5}
\caption{If $a=-1.1$ and $b=1$ then the largest three eigenvalues of the monodromy operator for equation \eqref{introex} are proven to be within the given discs.  Thus the parameter point $(a,b)$ is stable.  (Unit circle and circle $r=\delta=0.2$ also shown.)} \label{fig:introeigs}
\end{figure}

For instance, suppose $a=-1.1$ and $b=1$ in equation \eqref{introex}, a point indicated in figure \ref{fig:introchart} to be stable but near the stability boundary. Figure \ref{fig:introeigs} shows as solid dots the computed eigenvalues of the monodromy matrix.  For the largest three eigenvalues, discs are shown whose radii $r=.0434$ are \emph{proven} maximum distances by which the large eigenvalues of the monodromy operator might differ; see theorem \ref{thm:main}.  By ``large eigenvalues'' we mean that we have chosen to consider eigenvalues $\mu\in\CC$ such that $|\mu|>\delta =0.2$.

Because the largest eigenvalue has magnitude $|\mu_1|=.9369$ and thus $|\mu_1|+r<1$, figure \ref{fig:introeigs} represents a \emph{proof} that the parameters $(a,b)=(-1.1,1)$ are eigenvalue stable for equation \eqref{introex}.

\begin{figure}[ht]
\regfigure{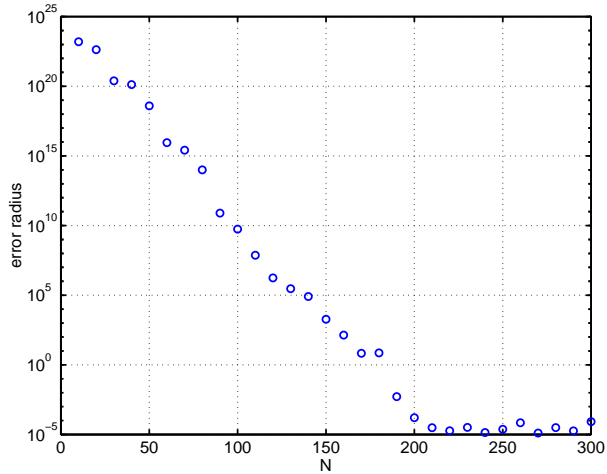}{2.5}
\caption{For DDE \eqref{introex} with $(a,b)=(-1.1,1)$: The monodromy matrix is an $(N+1)\times (N+1)$ matrix approximation to the $\infty\times \infty$ monodromy operator.  The radius found from theorem \ref{thm:main} of the error discs around the numerical eigenvalues of the monodromy matrix decays exponentially with $N$.} \label{fig:eradexpdecay}
\end{figure}

The numerical eigenvalues in figure \ref{fig:introeigs}, the dots, were generated using a $(N+1)\times (N+1)$ monodromy matrix with $N=184$.  It turns out that the radius of the discs decays roughly exponentially with $N$ as illustrated in figure \ref{fig:eradexpdecay}.  In this example, when $N\gtrapprox 220$ we find that floating-point errors, and the fact that the condition number of the matrix of eigenvectors of the monodromy matrix is roughly $10^6$, limits accuracy to roughly $10^{-5}$.

Subsection \ref{subsect:usethm} presents this scalar DDE example in detail.  Section \ref{sect:Mathieu} addresses a system of DDEs, a delayed Mathieu equation, with similar results.

This paper is somewhat more expository than is standard in the numerical analysis literature.  We start by recalling the basics of Chebyshev (spectral) collocation in section \ref{sect:interp}.  Then we prove new results for the first order scalar DDE case in sections \ref{sect:apost}, \ref{sect:DDEop} and \ref{sect:DDEapprox}.  We then generalize the same machinery to systems of DDEs (section \ref{sect:systems}) and apply the machinery to the example of a delayed, damped Mathieu equation (section \ref{sect:Mathieu}).  (In a less expository paper we would do the systems case immediately, with inevitable loss of clarity to the non-expert.)  Further generalizations and technical issues are covered in sections \ref{sect:polytechs} and \ref{sect:taunoteqperiod} and in the appendices. 

We now review some of the numerical analysis literature relating to stability of DDEs by spectral methods.  Spectral methods were perhaps first applied to DDE initial value problems in \cite{Bellen} and \cite{ITM}, though in the former work the emphasis is on ``$h$-refinement'' and in the latter the techniques are limited to constant-coefficient linear problems.  More recently, collocation using the Gauss-Legendre points, in particular, is applied to the problem of finding periodic solutions to nonlinear DDEs in \cite{ELHR} and the stability\footnote{We distinguish here between stability of the DDE and stability of the numerical method.} and convergence under $h$-refinement of the resulting method is addressed in \cite{EngelborghsDoedel}.

On the other hand, there exists some engineering literature addressing stability (as opposed to solutions of initial value problems) of ODEs and DDEs by spectral methods.  For ODE problems we specifically note \cite{LeeRenshaw} who use equally-spaced collocation points and \cite{SinhaWu,SSP} who use a Chebyshev Galerkin (or Galerkin-collocation, respectively) method.  In fact, such engineering applications led the current author and coworkers to DDE stability questions \cite{BMBAS}.

The current author knows of two \emph{numerical analysis} papers, namely \cite{LuzyaninaEngelborghs} and \cite{LER}, which address the application of spectral methods to the stability of linear DDEs and functional differential equations.  The former surveys collocation methods for approximating eigenvalues of monodromy operators (i.e.~Floquet multipliers) and gives numerical evidence for spectral convergence (see table 2 of \cite{LuzyaninaEngelborghs}).  In the latter work it is proven that the numerical method preserves a property (``RHP-stability'') which suffices for stability of the functional equation.  In neither case are estimates on the eigenvalues addressed.

This paper grows out of enjoyable work on DDEs in engineering applications initiated and sustained by Eric Butcher, and continued with students Victoria Averina, Tim Carlson, Haitao Ma, Praveen Nindujarla, Jake Stroh, and Ben White.

\section{Interpolation and spectral differentiation at Chebyshev points} \label{sect:interp}

We are interested in solving ODEs or DDEs on a fixed interval $t\in I$ and we suppose $I=[-1,1]$ as this choice is convenient relative to Chebyshev polynomial conventions.  All results are easily translatable to other intervals.  The functions we consider will generally be continuous and complex-valued; denote the space of such functions by $C(I)$.  (It will be most useful to consider complex-valued functions when considering stability problems as the eigenfunctions of the real monodromy operator (section \ref{sect:DDEop}) are generally complex.)  We use the norm $\|f\|_\infty = \max_{t\in I} |f(t)|$ for $f\in C(I)$.

\begin{defns}  Let $N\ge 1$.  Let $\Po_N$ be the space of at most degree $N$ polynomials with complex coefficients.  Let $\Co_N=\{t_0,\dots,t_N\}\subset I$ be the Chebyshev \emph{collocation points} (extreme points) $t_j=\cos(\pi j/N)$.  Note $t_0=1$, $t_N=-1$, and $t_{j+1}<t_j$.\end{defns}

As a subspace of $C(I)$ the dimension of $\Po_N$ is $N+1$.  For instance, the monomials $1,t,t^2,\dots,t^N$ are a basis of $\Po_N$, as are the Lagrange polynomials $L_j(t)=\prod_{k=0,\, k\ne j}^N (t-t_k)/(t_j-t_k)$ for the collocation points.

Polynomial interpolation using the collocation points is our first concern.  (See \cite{BurdenFaires} or \cite{Plato} for the basic theory of polynomial interpolation.)  Consider the following maps: 

\begin{defns}\hspace{1.0in}

\begin{itemize}
\item Let $\Ecal_N:C(I)\to \CC^{N+1}$ be the linear map which evaluates the input function at the $m$ collocation points $t_j$, so $(\Ecal_N f)_j = f(t_j)$ for $j=0,\dots,N$.
\item Let $E_N:\Po_N\to\CC^{N+1}$ be the restriction of $\Ecal_N$ to polynomials.  
\item Let $P_N:\CC^{N+1}\to \Po_N$ be the map which gives the unique polynomial of degree at most $N=m-1$ having the input values at the collocation points $t_j$.  Note $(P_N v)(t) = \sum_{j=0}^N v_j L_j(t)$, though we will have no practical need for this formula.
\item Let $I_N=P_N \Ecal_N: C(I)\to \Po_N$ be the interpolation map which constructs the unique polynomial of degree $N$ which agrees with the input function at the collocation points.\end{itemize}\end{defns}

Note that $E_N,P_N$ are linear isomorphisms, $P_N^{-1} = E_N$, and $I_N$ is a projection (i.e.~$I_N^2 = I_N$).

The quality of polynomial interpolation depends on the regularity of the function being interpolated.  Nonetheless it is possible to find \emph{analytic} functions $f$ on $I$ and (families of) interpolation points such that $f$ is not well-interpolated by polynomials \cite{Trefethen}.  Fortunately, Chebyshev collocation points are excellent for interpolation in the analytic case:

\begin{thm}[proven in \cite{Tadmor}, in particular]\label{thm:interp}  Suppose $f$ is analytic in an open, simply-connected region $R$ such that $[-1,1]\subset R \subset \CC$.  Then there exist constants $C>0$ and $\rho>1$ independent of $N$ such that if $p=I_N f\in\Po_N$ is the polynomial interpolant of $f$ at the $N+1$ Chebyshev collocation points $\Co_N$,
\begin{equation}\label{interperr}
   \|f-p\|_\infty \le C \rho^{-N}.
\end{equation}
Furthermore, suppose $E\subset R$ is an ellipse with foci $\pm 1$---see figure \ref{fig:ellipse} below.  If $S,s$ are the lengths of the semimajor, semiminor axes of $E$, respectively, then inequality \eqref{interperr} applies with $\rho=S+s$.  In fact, for $k=0,1,2,\dots$ there exists $C_k>0$ such that
    $$\|f^{(k)}-p^{(k)}\|_\infty \le C_k \rho^{-N},$$
for $\rho=S+s$, where $f^{(k)}$ is the $k$th derivative of $f$ on $I$. \end{thm}

\begin{figure}[ht]
\regfigure{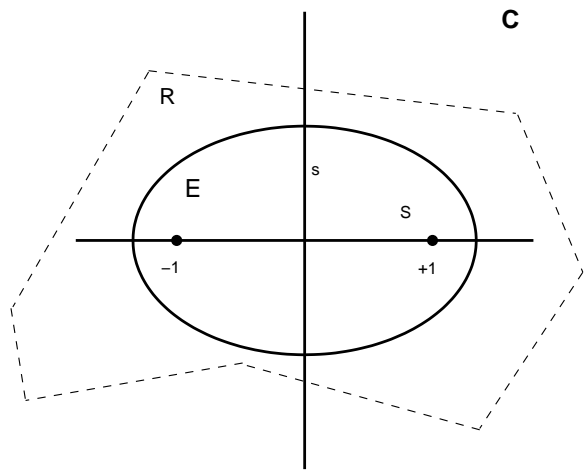}{2.5}
\caption{If $f$ is analytic in a region $R\subset \CC$ which contains a ``regularity ellipse'' $E$ with foci $\pm 1$ and semiaxes $S,s$ then Chebyshev interpolation converges at exponential rate $\rho^{-1}$ where $\rho=S+s$.  Note that $S^2=1+s^2$ for an ellipse with foci $\pm 1$.} \label{fig:ellipse}
\end{figure}

From polynomial interpolation one may construct spectral differentiation.  Let
    $$D_N=E_N \ddt{} P_N : \CC^{N+1} \to \CC^{N+1},$$
whose $(N+1)\times (N+1)$ matrix with respect to the standard basis for $\CC^{N+1}$ is called the \emph{Chebyshev spectral differentiation matrix}.  The matrix $D_N$ ``differentiates'' $v\in \CC^{N+1}$ as follows: construct a polynomial $p$ with values $v$ at the collocation points; differentiate it $p'=\frac{dp}{dt}$; and evaluate the result at the collocation points $(D_N v)_j=p'(t_j)$.  See \cite{Trefethen} for explicit formulas for the entries of $D_N$ and applications of $D_N$.  Note that all entries of $D_N$ are real and that $D_N$ is nilpotent.

Exponential decay of error as $N$ increases is known as \emph{spectral convergence}, so theorem \ref{thm:interp} says polynomial interpolation at the Chebyshev collocation points yields spectral convergence for analytic functions.  The goal in this note is to show spectral convergence for the solutions of ODEs and DDEs, and for certain eigenvalue problems, by the numerical methods described below.

We are interested in linear ODEs and DDEs which have highly regular (or piecewise-regular) coefficient functions.  In particular, because this assumption applies to the engineering applications which we have addressed and/or foresee \cite{BMBAS}, we suppose these coefficient functions are (piecewise) analytic.  In the DDE case we will additionally suppose that history functions are (piecewise) analytic. 

We end this section with a technical tool, namely, polynomial interpolation by discrete Chebyshev series.  Recall that $T_k(t)=\cos(k\arccos t)$ are the Chebyshev polynomials and that the collocation points $t_j=\cos(j\pi/N)$ are the extreme points of $T_N(t)$ since $T_N(t_j)=\cos(N j \pi/N)=(-1)^j$.  If $f\in C(I)$ then $p=I_N(f)$ is an $N$ degree polynomial which may be found by the following formulas
\begin{equation}\label{discreteCheb}
p(t)=\sum_{k=0}^N \tilde f_k T_k(t), \qquad \tilde f_k=\sum_{j=0}^N C_{kj} f(t_j), \qquad C_{kj} = \frac{2}{N \gamma_j \gamma_k} \cos\left(\frac{\pi j k}{N}\right),
\end{equation}
where $\gamma_j=2$ if $j=0,N$ and $\gamma_j=1$ otherwise.  These formulas may be implemented by a modification of the fast Fourier transform \cite{Trefethen}.  That is, the FFT may be used to compute $P_N$ and $E_N$ if desired.

\section{\emph{A posteriori} estimates for initial value problems}\label{sect:apost}

\subsection{An algorithm and error estimate}  \label{subsect:collocalg} Consider the scalar, linear ODE initial value problem
\begin{equation}\label{ODEnon}
\dot y(t)=a(t) y(t) + u(t), \qquad y(-1)=y_0.
\end{equation}
Needless to say, this problem has a well-known exact solution.  Let $\Phi_a(t)$ be the ODE fundamental solution solving $\dot \Phi_a(t) = a(t) \Phi_a(t)$ and $\Phi_a(-1)=1$.  In this scalar case $\Phi_a(t)=\exp\left(\int_{-1}^t a(s)\,ds\right)$, of course, but the existence of the fundamental solution $\Phi_a(t)$ generalizes to the systems case while this exponential formula does not (at least, directly \cite{Iserles}).  The exact solution to \eqref{ODEnon} is then
\begin{equation}\label{vopsoln}
y(t)=\Phi_a(t) \left[y_0 + \int_{-1}^t \Phi_a(s)^{-1} u(s)\,ds\right]
\end{equation}
by variation-of-parameters.  

Let $\hD_N$ be the $(N+1)\times (N+1)$ matrix modification of $D_N$ which satisfies $(\hD_N)_{ij}=(D_N)_{ij}$ for rows $i=1,\dots,N$ and all columns, but has last row
\begin{gather*}
    (\hD_N)_{N+1,k}=0, \quad 1\le k\le N, \qquad \text{and} \qquad (\hD_N)_{N+1,N+1}=1.
\end{gather*}
The product $\hD_N v$ represents both ``$\dot y(t)$'' for $-1<t\le 1$ and ``$y(-1)$'' if $v$ is the (Chebyshev) collocation approximation of $y(t)$.  For $a \in C(I)$ let $\hM_a$ be the $(N+1)\times (N+1)$ diagonal matrix with entries
    $$(\hM_a)_{ii} = \begin{cases} a(t_{i-1}), & 1\le i \le N,\\ 0, & i=N+1.\end{cases}$$
The product $\hM_a v$ represents the multiplication ``$a(t)y(t)$'' for $-1<t\le 1$ if $v$ is the collocation approximation of $y(t)$.  For $u\in C(I)$ and $y_0\in\CC$ let $\hu\in \CC^{N+1}$ have entries
\begin{equation}\label{uhatdefn}
\hu_i=\begin{cases} u(t_{i-1}), & 1\le i \le N, \\ y_0, & i=N+1.\end{cases}
\end{equation}
That is, $\hu_i=(\Ecal_N u)_i$ if $1\le i \le N$ while $\hu_{N+1}=y_0$.  The following is immediate from the construction of $\hD_N,\hM_a,$ and $\hu$.

\begin{lem}[the collocation algorithm]\label{lem:collocequiv}  Let $a,u\in C(I)$ and let $y_0\in \CC$.  The following two descriptions of the collocation approximation of the problem \eqref{ODEnon} are equivalent:
\renewcommand{\labelenumi}{(\textit{\roman{enumi}})}\begin{itemize}
\item $p\in \Po_N$ satisfies 
\begin{equation}\label{ODEmethpoly}
\dot p(t_j) = a(t_j) p(t_j) + u(t_j), \quad 0\le j \le N-1, \quad \text{ and } \quad p(-1)=y_0;
\end{equation}
\item $v\in\CC^{N+1}$ satisfies
\begin{equation}\label{ODEmeth}
\hD_N v = \hM_a v + \hu.
\end{equation}
\end{itemize}
The equivalence is $p=P_N v$ and $v=E_N p$.\end{lem}

We will call the numerical algorithm described by 
    $$y(t)\approx p(t) = P_N v= P_N \big(\hD_N-\hM_a\big)^{-1} \hu,$$
that is, by the above lemma, the \emph{collocation algorithm}.  (For a complete numerical algorithm we would need to specify techniques for constructing $D_N$ and for solving the linear equations which arise, but for now we assume these tasks are done exactly.  See, however, section \ref{sect:polytechs}.)

The algorithm generalizes gracefully to scalar DDEs (sections \ref{sect:DDEop} and \ref{sect:DDEapprox}) and to system of ODEs and DDEs.  It is well-proven in practice.

Our goal of showing spectral convergence for the collocation algorithm is, roughly speaking, the goal of showing that the collocation algorithm does as good a job of approximating the solution to the initial value problem as would polynomial interpolation on $y(t)$ in \eqref{vopsoln}.
 
What we can prove is an estimate on the maximum difference between $y(t)$ and $p(t)$ in terms of computable quantities.

\begin{thm}\label{thm:apostivp}  Suppose $y\in C(I)$ solves initial value problem \eqref{ODEnon} with $a\in C(I)$.  Suppose $p\in \Po_N$, $N\ge 1$, is calculated by the collocation algorithm (lemma \ref{lem:collocequiv}).  Let 
    $$C_a = \exp\left(\int_{-1}^1 \max\{\Re a(s),0\}\,ds\right)$$
and
    $$R_p=\dot p(-1) - a(-1) y_0 - u(-1).$$
Then the uniform error of $p$ as an approximation to $y$ satisfies
\begin{equation}\label{apostivp}
\|y-p\|_\infty \le 2 E_a \Big(\|ap-I_N(ap)\|_\infty + \|u-I_N(u)\|_\infty + |R_p|\Big).\end{equation}\end{thm}

Note that $|\Phi_a(t)| = \exp\left(\int_{-1}^t \Re a(s)\,ds\right) \le C_a$, so $C_a$ simply bounds the fundamental solution.
The quantity $R_p$ is a measure of the degree to which the approximating polynomial $p$ does not satisfy the differential equation at the initial time; it could be called the \emph{initial residual} for the method.  Actually computing the initial residual is easy, as $\dot p(-1) = (D_N v)_{N+1}$.

It would be nice to show that $\|ap-I_N(ap)\|_\infty$ and $R_p$ are (exponentially) small in $N$ and thus make this theorem into a (spectral) convergence result for the collocation algorithm.  For now, however, this theorem is ``\emph{a posteriori}'' because one must actually compute $p$ in order to evaluate the error estimate.

Computing the maximum norms $\|ap-I_N(ap)\|_\infty,\|u-I_N(u)\|_\infty$ is done by the techniques described in section \ref{sect:polytechs}.

In the special case where $a(t)=a_0$ is a constant, we can be more specific about the conclusion of this theorem.  If $\Re a_0>0$ then
\begin{equation}\label{specialone}
\|y-p\|_\infty \le 2 e^{2 \Re a_0} \|u-I_N(u)\|_\infty + c |R_p|
\end{equation}
where $c=(\pi(|a_0|+1)+4) \,e^{2\Re a_0}/(2\,N^2)$ and $R_p=\dot p(-1) -a_0 y_0 -u(-1)$.  If, $\Re a_0\le 0$, on the other hand, then
\begin{equation}\label{specialtwo}
\|y-p\|_\infty \le 2 \|u-I_N(u)\|_\infty + \min\left\{\tilde c,\frac{\pi}{2N}\right\} |R_p|,
\end{equation}
here $\tilde c=(\pi(|a_0|+1)+4) \,e^{-2\Re a_0}/(2\,N^2)$.  Note that the constant multiplying the initial residual actually decays with increasing $N$.  These inequalities are proven in appendix A.

The error estimate in theorem \ref{thm:apostivp} is absolute, not relative.  In practice one may choose $N$ so that $\|y-p\|_\infty \ll \|p\|_\infty$, established using the estimate, and then increase $N$, if needed, so that  $\|y-p\|_\infty/\|p\|_\infty < \delta$ for some relative tolerance $\delta$.  In the case that $\|p\|_\infty$ is itself small, absolute error $\|y-p\|_\infty$ is exactly what we would wish to control.

\subsection{Examples}\label{subsect:ivpexamples}  Before proving this theorem we wish to demonstrate its effectiveness in simple examples.  In particular, the behavior of different cases $a(t)=a_0$ (exponential growth when $\Re a_0\gg 0$, decay when $\Re a_0\ll 0$, and combined with oscillation when $|\Im a_0|$ is large) is captured by equations \eqref{specialone} and \eqref{specialtwo}.  We start with an easy case. 

\begin{figure}[ht]
\regfigure{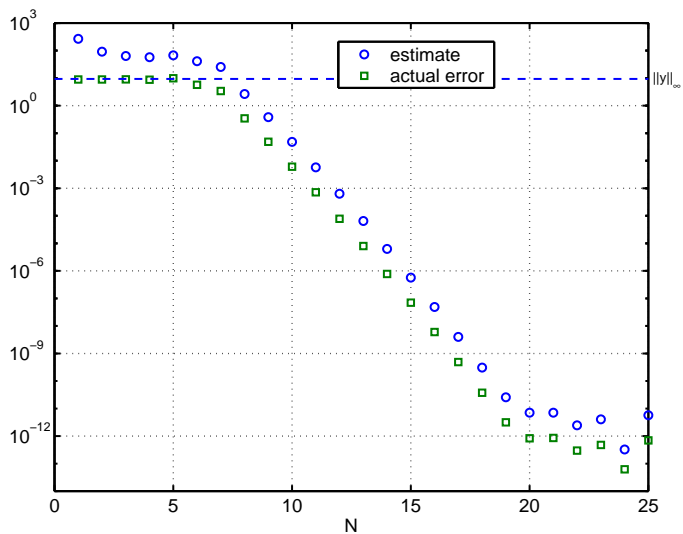}{2.8}
\caption{Actual error $\|y-p\|_\infty$ and the estimate \eqref{specialone} for the initial value problem in example \ref{exampleone}.} \label{fig:estvsactone}
\end{figure}

\begin{example}\label{exampleone}  Consider the following initial value problem
    $$\dot y(t) = 3 y(t) + t, \qquad  y(-1)=y_0,$$
with exact solution $y(t)=e^{3(t+1)} \left(y_0-\frac{2}{9}\right) - \frac{1}{3} \left(t +\frac{1}{3}\right)$.  We choose $y_0=0.2$ so that $y(t)$ takes values in the interval $[-9.5,0.2]$ for $t\in I$ and thus $y(t)$ is $O(1)$ on $I$.  Since $u(t)=t$ is a polynomial of degree one and $a(t)=3$ is constant, special case \eqref{specialone} gives the estimate 
    $$\|y-p\|_\infty \le 2 (\pi+1) e^6 N^{-2}  |\dot p(-1)+0.4|$$
for the maximum difference between the collocation polynomial $p\in\Po_N$ and the exact solution $y$.  We compare this estimate to the actual maximum error, computed by sampling the exact solution and approximating polynomial at $1000$ equally-spaced points in $I$.  See figure \ref{fig:estvsactone}.  Exponential decay of the actual error is clear.  The convergence stops roughly 3 orders of magnitude above machine precision $\eps_{\text{m}} \approx 2 \times 10^{-16}$.  Note that $\|y\|_\infty\approx 9.4$ so the relative error $\|y-p\|_\infty/\|y\|_\infty$ is roughly two orders of magnitude above $\eps_{\text{m}}$.  The estimate, however, is consistently only about about a factor of 10 too large, and this is quite acceptable because the convergence is so fast that the estimated value for $N$ to give a certain error is too big by only one in practice.\footnote{The estimate might be suspiciously good---it tracks the behavior of the actual error when both are comparable to machine precision---until one recalls the \emph{a posteriori} nature of the estimate.}\end{example}

\begin{figure}[ht]
\regfigure{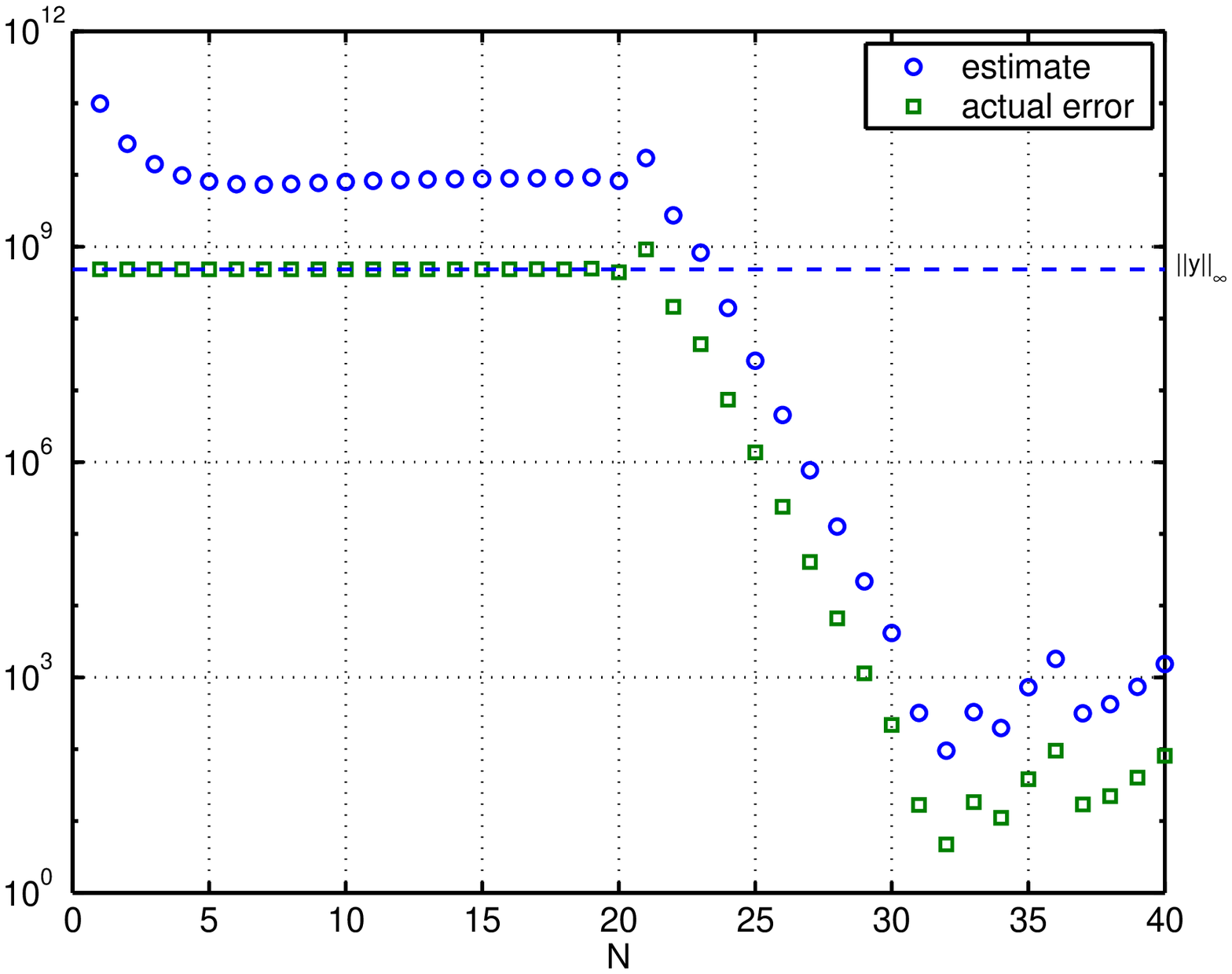}{2.3}\regfigure{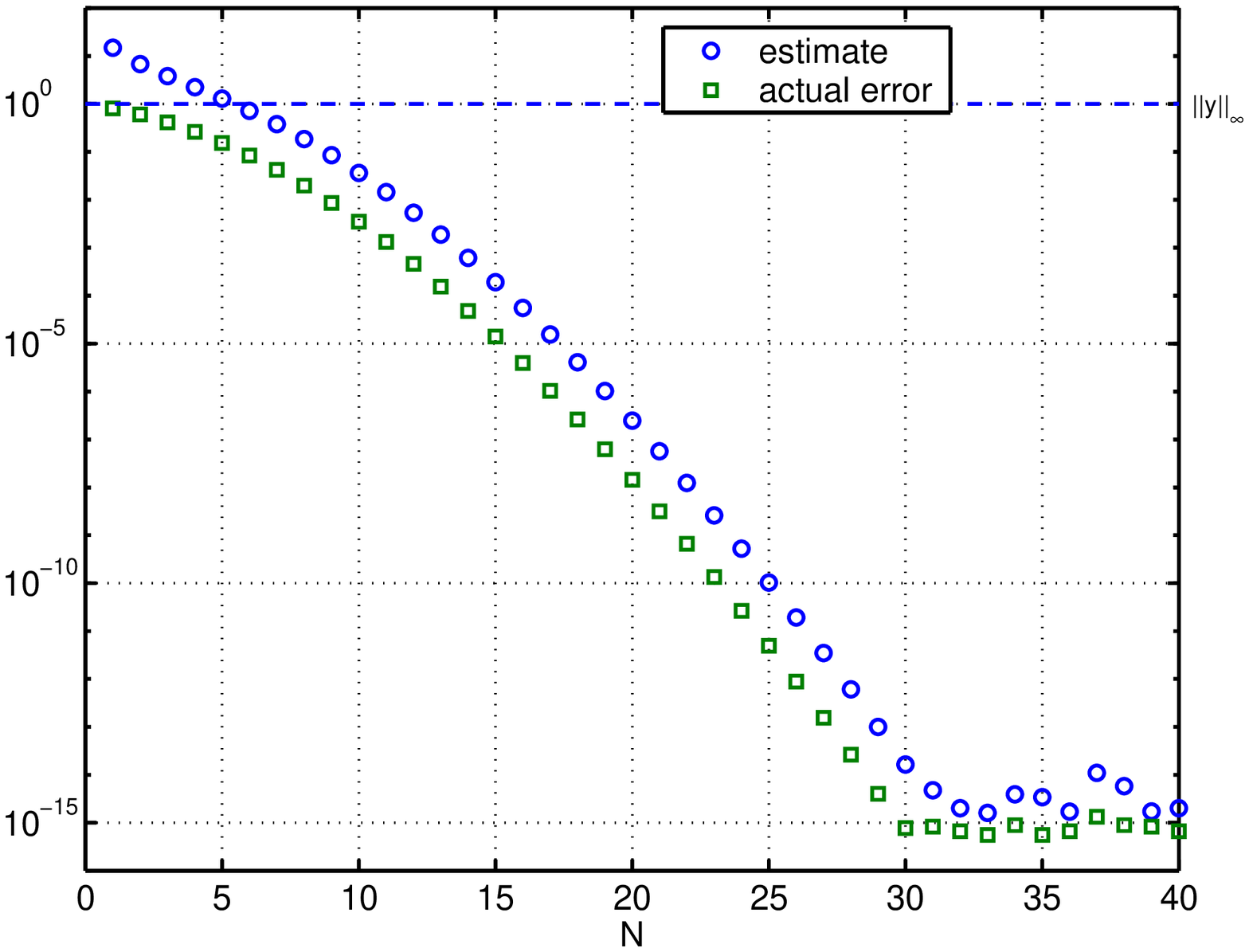}{2.3}
\caption{(a) Estimated and actual error for $\dot y = 10 y$.  (b) For $\dot y = -10 y$.} \label{fig:estvsacttwo}
\end{figure}

\begin{example}\label{exampletwo}  Cases with rapid exponential growth, decay, or oscillation---i.e.~stiffness---show that spectral methods have to work hard, too.  Let $u(t)=0$ and consider $\dot y(t) = a y(t)$, $y(-1)=1$, with exact solution $y(t)=e^{a(t+1)}$.  Figure \ref{fig:estvsacttwo} shows actual error and the estimate in cases $a=10$ and $a=-10$.  The behavior of the actual error is quite different in the two cases, and we use special cases \eqref{specialone} and \eqref{specialtwo}, respectively.  For $a=10$, note that the solution $y(t)=e^{10(t+1)}$ grows from $y(-1)=1$ to $y(1)=e^{20}\approx 5\times 10^8$, and thus the error is reasonable in a relative sense.  Also, convergence does not start until $N=20$ or so.  The significance of this number is that a polynomial of degree at least $20$ is required to approximate the rapid growth of $e^{10 t}$ on $I$.  For $a=-10$ we see different behavior, with immediate but slower convergence.  The explanation is that low order polynomials do exist which approximate $e^{-10t}$ in an absolute sense, but only polynomials of high degree succeed in being comparably small to $e^{-10 t}$ for $t\in I$ close to $+1$.  For both $a=\pm 10$ spectral convergence ends at about $N=30$, as non-truncation numerical errors (e.g.~errors in solving linear equations and rounding errors) become dominant. \end{example}

\begin{figure}[ht]
\regfigure{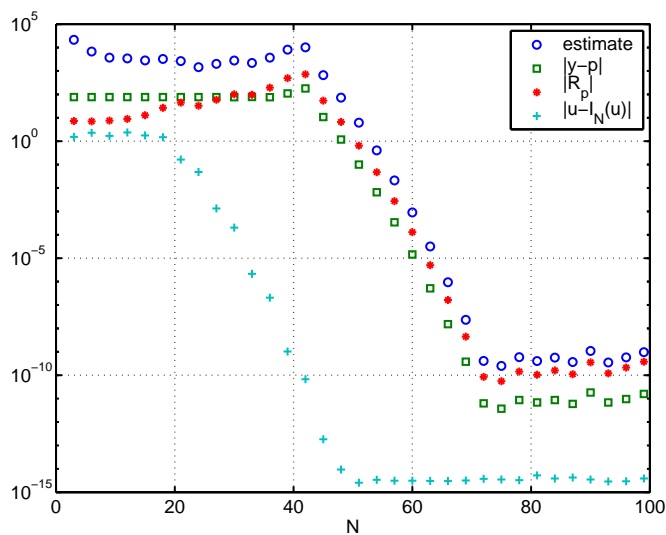}{2.8}
\caption{Estimated and actual error for $\dot y = (3+37i) y + \sin(20 t)$, $y(-1)=0.2$ for odd $N$ less than $100$.  Both terms on the right side of inequality \eqref{specialone} contribute to the estimate until $N\approx 50$.} \label{fig:estvsactthree}
\end{figure}

\begin{example}\label{examplethree}  The nonhomogeneous term  $u(t)$ is usually not a polynomial.  Suppose $a(t)=a_0=3+37 i$, $y(-1)=0.2$, and $u(t)=\sin(20 t)$.  The exact solution is found by variation-of-parameters.  We calculate the $\|u-I_N(u)\|_\infty$ term in theorem \ref{thm:apostivp} by evaluating $u$ and $I_N(u)$ at $1000$ equally spaced points.  Figure \ref{fig:estvsactthree} shows four quantities: the estimated error, the actual error, the initial time residual error $|\dot p(-1) - a y_0 - u(-1)|$, and the interpolation error $\|u-I_N(u)\|_\infty$.  The need for $N$ to be sufficiently large to ``handle'' the oscillations present in $\Phi_a(t)=e^{-a_0(t+1)}$ and $u(t)$ is clear. \end{example}

\begin{figure}[ht]
\regfigure{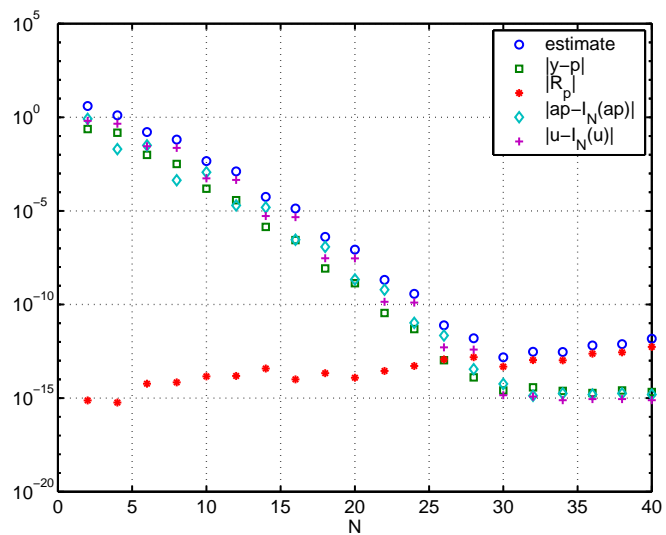}{2.8}
\caption{Estimated and actual error for $\dot y = 2 t \,y + t \sin(3 t^2)$, $y(-1)=1$.  Both interpolation error terms on the right side of inequality \eqref{apostivp} contribute to the estimate.} \label{fig:examplefour}
\end{figure}

\begin{example} Finally, the leading coefficient $a(t)$ need not be constant.  Suppose $a(t)=2t$, $u(t)=t\sin(3 t^2)$, $y(-1)=1$.  The exact solution can be found by variation of parameters because $\Phi_a(t)=e^{t^2-1}$ so that the relevant integral is $\int e^{1-s^2} s \sin(3 s^2)\,ds = e\int e^{-z} \sin(3 z)\,dz$ (under the substitution $z=s^2$).  Note $C_a=e$ is a conservative bound on $\Phi_a(t)$.  Figure \ref{fig:examplefour} shows that though the initial residual error happens to be small for all $N$, the interpolation errors $\|ap-I_N(ap)\|_\infty, \|u-I_N(u)\|_\infty$ are both active contributors to the error estimate and the actual error $\|y-p\|_\infty$.  These quantities decay exponentially till error and estimate are of close to machine precision $10^{-15}$.\end{example}

The collocation algorithm works well on initial value problems with or without the estimate given in theorem \ref{thm:apostivp}, including on ODE (and DDEs) systems and in cases with multiple delays and merely piecewise analytic coefficients.  We generally extract roughly double precision results---say, 8 to 14 correct figures for $\|a\|_\infty\lesssim 10$---in those problems for which we can compare to exact answers.  

A caveat: Theorem \ref{thm:apostivp} assumes that the spectral differentiation matrix $D_N$ is exactly computed and that linear system \eqref{ODEmeth} is exactly solved.  Computation with large $N$ (with $N\approx 100$ already fairly large) must be carefully handled to extend spectral convergence or even to maintain stability (see section \ref{sect:polytechs} and references \cite{BattlesTrefethen} and \cite{TrefethenTrummer}).

\subsection{Proofs}  
\begin{figure}[ht]
\regfigure{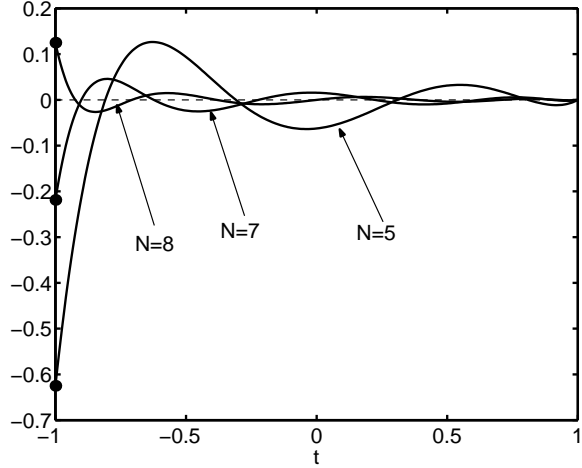}{2.5}
\caption{Polynomials $l_N(t)$ for three values of $N$.  It turns out that $l_N(t)$ has a closed-form expression in $\theta=\arccos t$.  The dots at $t=-1$ have values $(-1)^N N \, 2^{2-N}$.} \label{fig:lN}
\end{figure}

We now turn to proving theorem \ref{thm:apostivp}.  First we establish a lemma which is reminiscent of computations with the Dirichlet kernel in Fourier series \cite{Korner}, namely, we rewrite the monic polynomial $l_N(t)$ with roots $t_0,\dots,t_{N-1}$ as a closed-form trigonometric expression.  Figure \ref{fig:lN} shows graphs of the polynomials $l_N(t)$ for some values of $N$.

\begin{lem}\label{lem:lNtotrig}  For $N\ge 1$ let $l_N(t)=(t-t_0)\dots (t-t_{N-1})$.  For $t\in I$ and $t=\cos\theta$,
    $$l_N(t) = \frac{(\cos(\theta)-1) \sin(N\theta)}{2^{N-1} \sin(\theta)}, \qquad -1<t< 1,$$
and $l_N(-1)=(-1)^N N\, 2^{2-N}$, $\|l_N\|_\infty = N\, 2^{2-N}$, in particular.\end{lem}

\begin{proof}  We use discrete Chebyshev series (section \ref{sect:interp}) to write $l_N(t)$ as a sum of Chebyshev polynomials $T_k(t)$.  Noting that $l_N(t_j)=0$ for $j=0,\dots,N-1$,
\begin{align}
l_N(t) &= \frac{l_N(-1)}{2N} \left[1 - 2 T_1(t) + 2 T_2(t) - \dots + (-1)^{N-1} 2 T_{N-1}(t) + (-1)^N T_N(t)\right]. \label{lNasChebseries}
\end{align}

Now, since $l_N(t)$ is monic and because the coefficient of $t^N$ in $T_N(t)$ is $2^{N-1}$ (use the recursion relation $T_k=2 t T_{k-1} - T_{k-2}$), it follows that $1 = \frac{l_N(-1)}{2N} (-1)^N 2^{N-1}$ or $l_N(-1) = (-1)^N N\,2^{2-N}$ as claimed.  From equation \eqref{lNasChebseries} we have
    $$l_N(t) = \frac{(-1)^N}{2^{N-1}} \left[1 + 2 \sum_{k=1}^{N-1} (-1)^k \cos(k\theta) + (-1)^N \cos(N\theta)\right].$$
We produce a closed form expression for geometric sum,
\begin{align*}
1 + 2 \sum_{k=1}^{N-1} (-1)^k \cos(k\theta) &= \sum_{k=-(N-1)}^{N-1} \left(-e^{i\theta}\right)^k = (-1)^{N-1} \frac{\cos((N-1/2)\theta)} {\cos(\theta/2)},\end{align*}
and find, through trigonometric identities, that
\begin{align*}
l_N(t) &= \frac{(-1)^N}{2^{N-1}} \left[(-1)^{N-1} \frac{\cos((N-1/2)\theta)} {\cos(\theta/2)} + (-1)^N \cos(N\theta)\right] \\
    &= \frac{1}{2^{N-1}} \left[ \frac{-\sin(N\theta)\sin(\theta/2)} {\cos(\theta/2)}\right] = \frac{(\cos(\theta)-1)\sin(N\theta)}{2^{N-1} \sin(\theta)}.\end{align*} \end{proof} 

\begin{proof}[Proof of theorem \ref{thm:apostivp}]  Let $q=I_N(ap)$ and $w=I_N u$.  Note $q(t_j)=a(t_j)p(t_j)$ and $w(t_j)=u(t_j)$ for $j=0,\dots,N$, of course.

Let $r=\dot p - q - w$.  Note $r \in \Po_N$ and $r(t_j)=0$ for $j=0,\dots,N-1$ by lemma \ref{lem:collocequiv}.  Thus there is $\beta\in\CC$ such that $r(t) = \beta\, l_N(t)$.  By lemma \ref{lem:lNtotrig} and evaluating $r(t)$ with $t=-1$,
\begin{equation}\label{initialfact}
\beta l_N(-1) = \beta \frac{(-1)^N N}{2^{N-2}} = \dot p(-1) - a(-1) p(-1) - u(-1) = R_p.
\end{equation}

On the other hand, noting $\dot p = q + w + \beta\, l_N = ap + (q - ap + w +\beta\, l_N)$, variation-of parameters gives
\begin{equation}\label{voperreqn}
    p(t)=\Phi_a(t) \left[y_0 + \int_{-1}^t \Phi_a(s)^{-1} \left(q(s) - a(s) p(s) + w(s) + \beta\, l_N(s)\right)\,ds\right].
\end{equation}
Subtracting equations \eqref{vopsoln} and \eqref{voperreqn} yields
\begin{equation}\label{inttoestimate}
y(t)-p(t)=\int_{-1}^t \Phi_s^t\, \left[a(s)p(s)-q(s) + u(s)-w(s)\right] \,ds - \beta \int_{-1}^t \Phi_s^t\, l_N(s)\,ds.
\end{equation}
where $\Phi_s^t = \Phi_a(t)\Phi_a(s)^{-1}$.  Note that $|\Phi_s^t|\le C_a$ for $s\le t$.  Thus by lemma \ref{lem:lNtotrig},
    $$|y(t)-p(t)| \le 2 C_a \left(\|ap-q\|_\infty + \|u-w\|_\infty + |\beta| \frac{N}{2^{N-2}}\right).$$
Taking absolute values in equation \eqref{initialfact} we have proven the theorem.\end{proof}

Finally, a technical detail remains, needed when estimating eigenvalues of a DDE monodromy operator.

\begin{cor}\label{cor:apost}  Under the hypotheses of theorem \ref{thm:apostivp} we have the derivative estimate
    $$\|\dot y - \dot p\|_\infty \le (2 C_a \|a\|_\infty + 1) \Big(\|ap-I_N(ap)\|_\infty + \|u-I_N(u)\|_\infty + |R_p|\Big).$$
\end{cor}

\begin{proof}  Since $\dot y = ay+u$, $\dot p=q+w+\beta l_N$, using the notation of the above proof we have
    $$\dot y - \dot p = a(y-p) + (ap-q) + (u-w) + \beta l_N$$
so
    $$\|\dot y - \dot p\|_\infty \le \|a\|_\infty \|y-p\|_\infty + \|ap-q\|_\infty + \|u-w\|_\infty + |R_p|.$$
The result follows from theorem \ref{thm:apostivp}.\end{proof}

 \section{The monodromy operator for linear, periodic DDEs with integer delays} \label{sect:DDEop}

\subsection{The monodromy operator $U$ on $C(I)$} Consider the scalar linear DDE problem for $t\in I$:
\begin{equation}\label{DDEscalar}
    \dot x(t) = a(t) x(t) + b(t) x(t-2), \qquad x(s)=f(s+2) \text{ for } -3\le s \le -1.\end{equation}
We suppose $a,b,f\in C(I)$ and we extend $a,b$ periodically to $\RR$.  It is easy to solve problem \eqref{DDEscalar} by variation-of-parameters to find $x(t)$ corresponding to a given history function $f(t)$.  That is, for $t\in I$ we simply solve the ODE initial value problem $\dot x(t) = a(t) x(t) + u(t)$, $x(-1)=f(1)$ with $u(t)=b(t)f(t)$.  Suppose, furthermore, that $a(t),b(t)$ are periodic with period $T=2$.  The variation-of-parameters formula for \eqref{DDEscalar} can be translated to each interval of length $2$ to find $x(t)$ for all $t\ge -1$, and we regard the effect as an operator on $C(I)$.  

\begin{defn}  Define the \emph{DDE monodromy operator}\footnote{Or perhaps the \emph{delayed Floquet transition operator}.} $U:C(I)\to C(I)$ by
\begin{equation}\label{Udef}
(U f)(t) = \Phi_a(t) \left[f(1) + \int_{-1}^t \Phi_a(s)^{-1} b(s) f(s)\, ds\right].
\end{equation}
\end{defn}

In appendix B we show that $U$ is a compact operator and thus, in some sense, is as much like a matrix as any infinite rank operator.

The solution $x(t)$ can be found ``by steps'' with $x_n\in C(I)$, $n\ge 0$:
    $$x_{n+1} = U x_n, \qquad x_0=f.$$
The actual solution $x(t)$, $t\ge -1$, is built by ``concatenating'' the steps: $x(t)=x_n(t-2n-2)$ if $t\in[2n-3,2n-1]$.

DDE \eqref{DDEscalar} is therefor stable if and only if $U^k f \to 0$ as $k\to\infty$ for all $f\in C(I)$.  Equivalently, if $\rho(U)=\max_{\lam\in\sigma(U)} |\lam|$ is the spectral radius of $U$ \cite{Lax}, $\rho(U)<1$ if and only if $U^k f \to 0$ as $k\to\infty$ for all $f\in C(I)$. This is our stability criterion, an eigenvalue criterion though there exist more sophisticated measures \cite{TrefethenPO}.

We are thus motivated to consider eigenfunctions of $U$.  Suppose $Ux=\lam x$ for $\lam \ne 0$.  Then $y=\lam x$ solves the ODE $\dot y = a y + b x$, that is, 
\begin{equation}\label{eigfcnODE}
\dot x = \left(a + \frac{1}{\lam} b\right) x.
\end{equation}
In this scalar case the solution is
\begin{equation}\label{eigform}
x(t)=x(-1)\exp\left(\int_{-1}^t a(s)+\frac{1}{\lam}b(s)\,ds\right).
\end{equation}
Of course, such a formula does not help us find $\lam$ or $x$ immediately, but we will use it to address the quality of polynomial interpolation of $x$ in the next section.  Note, however, that since $(Uf)(-1)=f(1)$, $\lam = \exp(\tilde a + \tilde b/\lam)$ where $\tilde a = \int_{-1}^1 a(s)\,ds$ and $\tilde b = \int_{-1}^1 b(s)\,ds$.  This allows determination of eigenvalues by complex variable methods.  \label{easycasepage}  Such methods do not generalize to systems, however, because the fundamental solution $\Phi_{a,b,\lam}(t)$ to equation \eqref{eigfcnODE} is not generally found by straightforward integration \cite{Iserles}.

\subsection{A useful Hilbert space ${H^1}$}\label{subsect:Hone}  Theorem \ref{thm:BFH} below, which will allow us to compare the eigenvalues of $U$ to those of $U_N$, requires operators acting on a Hilbert space.  The obvious space is the one for which the Chebyshev polynomials form an orthogonal basis.

\begin{defn}  Let $L^2=L^2_T(I)$ be the Hilbert space of complex-valued (measurable) functions $f$ on $I$ such that $\int_{-1}^1 |f(t)|^2 (1-t^2)^{-1/2}\,dt < \infty$.  (Where ``T'' is for ``Chebyshev''.)  Let 
    $$\ip{f}{g}_{L^2} = \int_{-1}^1 \overline{f(t)}\, g(t)\, \left(1-t^2\right)^{-1/2}\,dt, \qquad \|f\|_{L^2}^2 = \ip{f}{f}_{L^2}.$$   \end{defn}

\begin{defn}  Define the (normalized) Chebyshev polynomials $\hT_0(t)=(1/\sqrt{\pi}) T_0(t)=1/\sqrt{\pi}$, $\hT_k(t)=(\sqrt{2/\pi}) T_k(t) = \sqrt{2/\pi} \cos(k\arccos t)$.  The set $\{\hT_k\}_{k=0}^\infty$ is an orthonormal (``ON'') basis of $L^2$.  For $f\in L^2$ let $\hf_k=\ip{\hT_k}{f}_{L^2}$, the \emph{Chebyshev expansion coefficients} of $f$, so $f(t)=\sum_{k=0}^\infty \hf_k \hT_k(t)$ (with convergence in $L^2$) and $\|f\|_{L^2}^2 = \sum_{k=0}^\infty |\hf_k|^2$.\end{defn}

We might hope that $U$ is a nice operator on $L^2$ with matrix entries we can calculate.  Unfortunately, $U$ is not even \emph{bounded} on $L^2$.  This is because formula \eqref{Udef} defining $U$ refers to the pointwise values of $f$, while elements of $L^2$ are actually equivalence classes which can take any values on given sets of measure zero.

\begin{exer}  Assume $b(t)\equiv 0$.  Construct a sequence of functions $f_n\in C(I)$ such that $f_n(1)=1$ and $\|f_n\|_{L^2}\to 0$ but $\|Uf_n\|_{L^2} \not\to 0$.\end{exer}

Thus we need to refine our choice of Hilbert space.

\begin{defn}  We define $H^1=H_T^1(I)$ in terms of the expansion coefficients,
    $$H^1 \equiv \left\{f\in L^2 \Big| \sum_{k=0}^\infty (1+k)^2 |\hf_k|^2 <\infty\right\} \subset L^2.$$
For $f,g \in H^1$ we define the inner product and norm
    $$\ip{f}{g}_{H^1} = \sum_{k=0}^\infty (1+k)^2 \overline{\hf_k} \,\hg_k, \qquad \|f\|_{H^1}^2 = \ip{f}{f}_{H^1}.$$
Let 
    $$\tilde T_k(t) = (1+k)^{-1} \hT_k(t) = \begin{cases} \sqrt{1/\pi}\, T_0(t), & k=0, \\ \sqrt{2/\pi}\, (1+k)^{-1} T_k(t), & k\ge 1,\end{cases}$$
the $H^1$-normalized Chebyshev polynomials.  The set $\{\tilde T_k\}_{k=0}^\infty$ is an ON basis of $H^1$.\end{defn}

The Hilbert space $H^1$ is introduced in \cite{Tadmor}.  It is a Sobolev space  but it is not equivalent to $W_T^1 = \left\{f\in L^2 \big|\|f\|_{L^2} + \|\dot f\|_{L^2}<\infty\right\}$.  Elements of $H^1$ are, more regular that those of $L^2$.  In fact, $H^1 \subset C(I)$ and pointwise evaluation is bounded.
 
\begin{lem}\label{lem:pointwiseeval}  If $f\in{H^1}$ then $f$ is continuous (has a continuous representative $\tilde f$) and $|\tilde f(t)|\le 0.9062 \|f\|_{H^1}$.
In particular,
\begin{equation}\label{deltaone}
\delta_1 f \equiv \sum_{k=0}^\infty \hf_k \hT_k(1) = \sum_{k=0}^\infty \int_{-1}^1 \hT_k(1) \hT_k(t) (1-t^2)^{-1/2} f(t)\,dt
\end{equation}
is a bounded linear functional on $f$ with $\|\delta_1\| \le 0.9062$ and $\delta_1 f = \tilde f(1)$.\end{lem}

\begin{proof}  The sum defining $\delta_1 f$ converges absolutely, and, in fact, $\sum_{k=0}^\infty \hf_k \hT_k(t)$ converges uniformly for $f\in{H^1}$:
\begin{align*}
\sqrt{\pi/2} \sum_{k=1}^\infty |\hf_k| |\hT_k(t)| &\le \sum_{k=1}^\infty |\hf_k| =  \sum_{k=1}^\infty \frac{1}{1+k} \cdot (1+k) |\hf_k| \\
    &\le \left(\sum_{k=1}^\infty (1+k)^{-2}\right)^{1/2} \left(\sum_{k=1}^\infty (1+k)^2 |\hf_k|^2\right)^{1/2}<\infty
\end{align*}
using Cauchy-Schwarz.  It follows that $f$ is continuous because (in addition) $\sum \hf_k \hT_k(t) = f$ a.e.   Note $\sum_{k=1}^\infty (1+k)^{-2}=\pi^2/6 - 1$.  Thus
\begin{align*}
|f(t)| &\le \sum_{k=0}^\infty |\hf_k| |\hT_k(t)| = \sqrt{1/\pi} |\hf_0| + C \left(\sum_{k=1}^\infty (1+k)^2 |\hf_k|^2\right)^{1/2} \\
    &\le \left((2/\pi)|\hf_0|^2 + 2 C^2 \sum_{k=1}^\infty (1+k)^2 |\hf_k|^2\right)^{1/2} \le \sqrt{2}\,C \|f\|_{H^1}
\end{align*}
where $C^2=\pi/3 - 2/\pi$ and using $a+b\le \sqrt{2 a^2 + 2 b^2}$.  Note $\sqrt{2}\,C \le 0.9062$. \end{proof}

Equation \eqref{deltaone} is meant to suggest $\text{``}\delta_1(t) = \sum_{k=0}^\infty \hT_k(1) \hT_k(t) (1-t^2)^{-1/2}\text{''}$ as a ``delta function'', but this sum does not converge.

In using the Hilbert space ${H^1}$ we will need to know that if $f$ is sufficiently regular on $I$ then $f\in{H^1}$.  We give a criteria via the Fourier series of $f(\cos\theta)$.

\begin{defn}  Let $L_F^2$ be the space of measurable functions $g$ on $[-\pi,\pi]$ such that $\int_{-\pi}^\pi |g(\theta)|^2\,d\theta <\infty$.  (Here ``F'' is for ``Fourier''.)  Let $\|g\|_F^2 = \int_{-\pi}^\pi |g(\theta)|^2\,d\theta$.\end{defn}

\begin{lem}\label{lem:ConethenHcal}  Suppose $f\in C^1(I)$.  The even function $\uf(\theta)=f(\cos\theta)$ is in $C_{per}^1[-\pi,\pi]$, that is, $\uf$ can be periodically extended with period $2\pi$ to be $C^1$ on $\RR$.  For $k\in\ZZ$ let $\hat{\uf}(k)$ be the $k$th Fourier coefficient of $\uf$, that is, $\hat{\uf}(k) = (2\pi)^{-1/2} \int_{-\pi}^\pi e^{-ik\theta} \uf(\theta)\,d\theta$.  Then $\hf_k = \hat{\uf}(k) = \hat{\uf}(-k)$ for $k>0$ and $\hf_0 = 2^{-1/2} \hat{\uf}(0)$.  

It follows that
\begin{align}
    \|f\|_{L^2}^2 &= \frac{1}{2} \sum_{k=-\infty}^\infty |\hat{\uf}(k)|^2 = \frac{1}{2} \|\uf\|_F^2, \text{ and } \label{FourLtwo} \\
    \|f\|_{{H^1}}^2 &= \frac{1}{2} \sum_{k=-\infty}^\infty (1+k)^2|\hat{\uf}(k)|^2 \le \|\uf\|_F^2 + \|\uf'\|_F^2.\label{FourHcal}
\end{align}
Thus $C^1(I) \subset {H^1} \subset C(I)$ and
\begin{equation}\label{FourinfHcal}
\|f\|_{H^1}^2 \le 2\pi \|f\|_\infty^2 + 2\pi \|\dot f\|_\infty^2
\end{equation}
if $f\in C^1(I)$.\end{lem}

\begin{proof}  The first assertion is elementary calculus.  Now
    $$\hf_k = \ip{\hT_k}{f}_{L^2} = \sqrt{\frac{2}{\pi}} \int_0^\pi \cos(k\theta) \uf(\theta)\,d\theta = \frac{1}{2} \sqrt{\frac{2}{\pi}} \int_{-\pi}^\pi \cos(k\theta) \uf(\theta)\,d\theta$$
if $k>0$ because $\uf$ is even.  Because $\uf$ is even it also follows that $\hat{\uf}(k)$ is real and thus $\hf_k=\hat{\uf}(k)$.  Similarly for $\hf_0$.  Using Parseval's equality, equation \eqref{FourLtwo} follows from
    $$\|\uf\|_F^2 = \sum_{k=-\infty}^\infty |\hat{\uf}(k)|^2 = |\hat{\uf}(0)|^2 + 2 \sum_{k=1}^\infty |\hat{\uf}(k)|^2 = 2 \sum_{k=0}^\infty |\hf_k|^2.$$
Next, if $g\in C_{per}^1[-\pi,\pi]$ then $\|g'\|_F^2 = \sum_{k=-\infty}^\infty k^2 |\hat g(k)|^2<\infty$, and thus \eqref{FourHcal}.  Finally, since $\uf'(\theta) = -\dot f(\cos\theta) \sin\theta$, inequality \eqref{FourHcal} gives \eqref{FourinfHcal}:
    $$\|f\|_{H^1}^2 \le \int_{-\pi}^\pi |f(\cos\theta)|^2\,d\theta + \int_{-\pi}^\pi |\dot f(\cos\theta)|^2 \sin^2\theta\,d\theta \le 2 \pi \|f\|_\infty^2 + 2\pi \|\dot f\|_\infty^2.$$\end{proof}

\subsection{The monodromy operator $\hU$ on ${H^1}$} \label{subsect:UonHilbert} We can now define the monodromy operator acting on ${H^1}$.  We then find its norm and estimate the accuracy of polynomial interpolation of its eigenfunctions, crucial steps in using theorem \ref{thm:BFH} below.

\begin{defn} Suppose $a,b\in C(I)$ in DDE \eqref{DDEscalar}.  Define $\hU\in\Lop({H^1})$ by
\begin{equation}\label{hUdefn}
    \left(\hU f\right)(t) = \Phi_a(t)\left[f(1) + \int_{-1}^t \Phi_a(s)^{-1} b(s)f(s)\,ds\right]
\end{equation}
for $f\in {H^1}$.\end{defn}

The eigenvectors of $U$ and $\hU$ are identical.  In fact, let $i_{H^1}:H^1 \hookrightarrow C(I)$ be the injection and note:

\begin{lem} $U,\hU$ are related by $\hU=U\circ i_{H^1}$ on ${H^1}$ and furthermore $Uf=\lam f$ for $f\in C(I)$ and $\lam\ne 0$ if and only if $f\in{H^1}$ and $\hU f =\lam f$.\end{lem}

\begin{proof}  Note that if $f\in C(I)$ then $Uf\in C^1(I)\subset {H^1}$ by definition of $U$ (formula \eqref{Udef}) and lemma \ref{lem:ConethenHcal}.  But if $Uf=\lam f$ it follows that $f\in{H^1}$ and thus $\hU f = U f =\lam f$.  Conversely, if $f\in{H^1}$ then $f\in C(I)$ by lemma \ref{lem:pointwiseeval} and thus $Uf=\hU f$.\end{proof}

In appendix B we see that the compactness of $\hU$ follows from that of $U$.  On the other hand we have an \emph{a priori} estimate of the norm of $\hU$.

\begin{lem}\label{lem:Uhatnorm}  If $C_a = \exp\left(\int_{-1}^1 \max\{\Re a(s),0\}\,ds\right)$ then 
    $$\|\hU\| \le c_0 + c_1 C_a + c_2 C_a^2$$
as an operator on ${H^1}$, where $c_0=\|b\|_\infty$, $c_1=2.3(1+\|a\|_\infty)+\pi\|b\|_\infty$, and $c_2=\pi\sqrt{2}\|a\|_\infty\|b\|_\infty$.\end{lem}

\begin{proof}  Clearly $\|\hU f\|_{H^1} \le \|\Phi_a\|_{H^1} |\delta_1 f| + \|\Phi_a g \|_{H^1}$ where $g(t)=\int_{-1}^t \Phi_a(s)^{-1} b(s) f(s)\,ds$.  Note $|\Phi_a(t)|\le C_a$ and $|\dot\Phi_a(t)|\le C_a \|a\|_\infty$.  (Also $|\Phi_a(t) \Phi_a(s)^{-1}|\le C_a$ if $s\le t$; this is useful in \eqref{maxparts} below.)  Thus $\|\Phi_a\|_{H^1} \le \sqrt{2\pi} C_a \left(1+\|a\|_\infty^2\right)^{1/2}$ by \eqref{FourinfHcal}.  It follows that 
    $$\|\Phi_a\|_{H^1} |\delta_1 f| \le 2.2715 (1+\|a\|_\infty) C_a \|f\|_{H^1}$$
by lemma \ref{lem:pointwiseeval} and using $(a^2+b^2)^{1/2}\le a+b$ for $a,b\ge 0$.

Now note
\begin{equation}\label{maxparts}
\max\{|\Phi_a(t)g(t)|,|g(t)|\} \le C_a \|b\|_\infty \int_{-1}^1 |f(s)|\,ds \le \sqrt{\pi/2}\, C_a \|b\|_\infty \|f\|_{L^2}
\end{equation}
using Cauchy-Schwarz.  Again from \eqref{FourinfHcal}, noting $\Phi_a(t)\dot g(t) = b(t) f(t)$,
\begin{align*}
\|\Phi_a g\|_{H^1}^2 &\le 2\pi \|\Phi_a g\|_\infty^2 + 2\pi \|\dot\Phi_a g + \Phi_a \dot g\|_\infty^2 \\
    &\le \pi^2 C_a^2 \|b\|_\infty^2 \|f\|_{L^2}^2 + 2\pi \left(C_a \|a\|_\infty \sqrt{\pi/2}\,\|b\|_\infty \|f\|_{L^2} + \|b\|_\infty \|f\|_\infty\right)^2 \\
    &\le \pi^2 C_a^2 \|b\|_\infty^2 \|f\|_{H^1}^2 + 4\pi \left((\pi/2) C_a^4 \|a\|_\infty^2 \|b\|_\infty^2 \|f\|_{H^1}^2 + (0.9062)^2 \|b\|_\infty^2 \|f\|_{H^1}^2\right) \\
    &= \|b\|_\infty^2 \left((0.9062)^2 + \pi^2 C_a^2 + 2\pi^2 \|a\|_\infty^2 C_a^4\right) \|f\|_{H^1}^2,
\end{align*}
using $(a+b)^2 \le 2 a^2 + 2 b^2$, the triviality $\|f\|_{L^2}\le \|f\|_{H^1}$, and lemma \ref{lem:pointwiseeval}.  Therefore
    $$\|\hU\|_{H^1} \le 2.2715 (1+\|a\|_\infty) C_a + \|b\|_\infty \left((0.9062)^2 + \pi^2 C_a^2 + 2\pi^2 \|a\|_\infty^2 C_a^4\right)^{1/2}$$
and the claim follows using $(a^2 + b^2 + c^2)^{1/2} \le a+b+c$ for $a,b,c\ge 0$.\end{proof}

We now turn to the regularity of the eigenfunctions of $\hU$.  Recall that if $\hU x = \lam x$, $\lam \ne 0$ then formula \eqref{eigform} applies.

\begin{lem}\label{lem:eiginterp}  Suppose $k\ge 1$.  Suppose $a(z),b(z)$ are analytic functions on a common (closed) ellipse $E\subset \CC$ with foci $\pm 1$ and sum of semiaxes $e^\eta=S+s>1$ (see theorem \ref{thm:interp}).  Let $\|a\|_{\infty,E}=\max_{z\in E} |a(z)|$ and similarly define $\|b\|_{\infty,E}$.  Suppose $\hU x=\lam x$ for $x\in {H^1}$ and $\lam\ne 0$, and suppose $\|x\|_{{H^1}} = 1$.  Let $A_E=\max_{z\in E} \left|\int_{-1}^z a(\zeta)\,d\zeta\right|$, $B_E=\max_{z\in E} \left|\int_{-1}^z b(\zeta)\,d\zeta\right|$.  For $p=I_{k+1} x\in \Po_k$,
\begin{equation}\label{eigfcninterp}
\|x-p\|_{H^1} \le 8 (\sinh\eta)^{-1} \exp\left(A_E + B_E/|\lam|\right)\, k \, e^{-k\eta}. 
\end{equation}\end{lem}

\begin{proof}  First, formula \eqref{eigform} determines an analytic continuation $x(z)$ of $x(t)$ defined on $I$ to $z\in E$.   Apply theorem 4.4 and, specifically, inequality (4.16) in \cite{Tadmor} to $x(z)$, noting that
    $$\max_{z\in E} \left|x(z)\right| = |x(-1)| \exp\left(\max_{z\in E}\left|\int_{-1}^z a(s) + \frac{1}{\lam} b(s)\,ds\right|\right) \le |x(-1)| \exp\left(A_E + B_E/|\lam|\right),$$
and conclude $\|x-p\|_{H^1} \le |x(-1)| 8 (\sinh\eta)^{-1} \exp\left(A_E + B_E/|\lam|\right) k e^{-k\eta}$.  Also  $|x(-1)|\le 0.9062$ by lemma \ref{lem:pointwiseeval}.  Estimate \eqref{eigfcninterp} follows.\end{proof}

 \section{Approximating the DDE monodromy operator and its eigenvalues} \label{sect:DDEapprox}

\subsection{Approximation of $U$}  We propose a Chebyshev spectral collocation method for the linear DDE initial value problem \eqref{DDEscalar} based upon the ODE method in section \ref{sect:apost}.   Namely, the approximation $p\in \Po_N$ to $x\in C(I)$ which solves \eqref{DDEscalar} is defined by
\begin{equation}\label{DDEmethpoly}
    \hD_N v = \hM_a v + \hu, \quad \text{ where } \quad \hu_j = \begin{cases} b(t_j) f(t_j), & j=0,1,\dots,N-1, \\ f(1), & j=N,\end{cases}
\end{equation}
where $v\in\CC^{N+1}$ and $\hD_N,\hM_a$ are $(N+1)\times (N+1)$ matrices as before.  As in lemma \ref{lem:collocequiv}, $v$ is equivalent to $p$: $v=E_N p$, $p=P_N v$.  

In fact, there is an $(N+1)\times (N+1)$ matrix approximation to the delayed term ``$b(t) x(t-2)$'' in equation \eqref{DDEscalar} defined by
    $$(\hM_b w)_i = \begin{cases} b(t_{i-1}) w_i, & 1\le i \le N, \\ w_1, & i=N+1,\end{cases}$$
where $w\in\CC^{N+1}$.  Compare equation \eqref{uhatdefn} defining $\hu$.

The collocation method now approximates $U$ itself and not merely solutions to particular initial value problems.  Let
\begin{equation}\label{Umdef}
U_N = \left(\hD_N-\hM_a\right)^{-1} \hM_b,
\end{equation}
a linear operator on $\CC^{N+1}$.  (Existence of the inverse in \eqref{Umdef} is not automatic, but one can show that for continuous $a(t)$ it exists for all $N$ sufficiently large \cite{TrefethenTrummer}.) 

The finite-dimensional linear map $U_N$ is an approximation to $\hU$ on ${H^1}$ as we will show in an \emph{a posteriori} sense.  To be more precise we need to redefine $U_N$ to act on ${H^1}$.  Recall that ${H^1}\subset C(I)$.

\begin{defn}  Let $\Pi_N$ be the orthogonal projection on ${H^1}$ which has range $\Po_N$, the degree $N$ polynomials.  Define $\hU_N\in\Lop({H^1})$ by
    $$\hU_N \equiv P_N U_N E_N \Pi_N.$$
\end{defn}

Clearly $\operatorname{rank}(\hU_N)\le N+1$.  Noting that $\operatorname{range}(\hU_N)\subset \Po_N \subset {H^1}$, it follows that $\hU_N p = \lam p$ if and only if $p\in\Po_N$ and $U_N v = \lam v$ where $v=E_N p$.  We will show in theorem \ref{thm:main} that because (in an \emph{a posteriori} sense) $\hU_N f \approx \hU f$ on a finite-dimensional subspace of $H^1$, the large eigenvalues of $\hU_N$ approximate those of $\hU$.

\subsection{Eigenvalue perturbation for diagonalizable operators on Hilbert spaces}  \label{subsect:eigpert}  The perturbation result we prove and use, theorem \ref{thm:BFH}, generalizes the well-known Bauer-Fike theorem for matrices \cite{BauerFike}, discussed in appendix C.    Theorem \ref{thm:BFH} is probably not new and is fairly close to the core idea of the eigenvalue condition number introduced in \cite{Wilkinson}.  It is also close to proposition 1.15 in \cite{Chatelin}.

Recall that a separable infinite-dimensional Hilbert space is isometrically isomorphic to $l^2$ (see, e.g. \cite{Lax} for this basic fact) where $l^2$ denotes the space of sequences $a=(a_1,a_2,\dots)$, $a_j\in\CC$, such that $\sum |a_n|^2 < \infty$.  Denote the standard basis elements of $l^2$ by $\delta_j$: $(\delta_j)_n=1$ if $n=j$ and $(\delta_j)_n=0$ otherwise.  Note that $\Lambda\in \Lop(l^2)$ is \emph{diagonal in the standard basis} if for each $j$, $\Lambda \delta_j = \lam_j \delta_j$ for some $\lam_j\in\CC$.

\begin{thm} \label{thm:BFH} Let $\Xcal$ be a separable infinite-dimensional complex Hilbert space (or, in the finite-dimensional case, suppose $\Xcal\cong \CC^d$).  Suppose $A\in\Lop(\Xcal)$, a bounded operator on $\Xcal$, is diagonalizable in the sense that there is a linear, bounded map $V:l^2 \to \Xcal$ (respectively, $V:\CC^d \to \Xcal$) with bounded inverse and a diagonal (in the standard basis) operator $\Lambda\in\Lop(l^2)$ (respectively, $\Lambda\in\CC^{d\times d}$) such that $A=V \Lambda V^{-1}$.  If $B\in\Lop(\Xcal)$ and if $Bx=\mu x$ for $\mu\in\CC$ and $\|x\|=1$ then 
\begin{equation}\label{BFestHilb}
\min_{\lam\in\sigma(A)} |\mu-\lam| \le \|V\| \|V^{-1}\| \|(B-A)x\|.
\end{equation}
\end{thm}

The hypotheses of theorem \ref{thm:BFH} imply that $\sigma(A)$ is the closure of the set of diagonal entries of $\Lambda$: $\sigma(A) = \overline{\{\lam_j\}}$.  Recall that $\sigma(A)\subset\CC$ is compact for any bounded operator $A\in\Lop(\Xcal)$, and thus the ``$\min$'' in estimate \eqref{BFestHilb} is appropriate.  On the other hand, the hypotheses do not imply that $A$ is a compact operator.

Theorem \ref{thm:BFH} can be generalized to Banach spaces which are isomorphic to sequence spaces (i.e.~$l^p$ spaces).  More specifically, we need to hypothesize an operator $\Lambda$, similar to $A$, for which $\|(\Lambda-\mu I)^{-1}\| \le \sup_{\lambda_i} |\lambda_i-\mu|^{-1}$ where $\{\lambda_i\}$ is a dense subset of $\sigma(A)$.  It is clear that assuming $A$ similar to a diagonal operator on $l^p$ suffices, for instance.  On the other hand, if $\Lambda$ is merely normal \cite{Lax} on a Hilbert space then the just-mentioned spectral property will hold.

Each norm appearing on the right side of estimate \eqref{BFestHilb} is different.  Recall that if $\Vcal_1$, $\Vcal_2$ are normed vector spaces with norms $\|\cdot\|_1$, $\|\cdot\|_2$ and if $L:\Vcal_1\to \Vcal_2$ is linear and bounded then $\|L\|_{\Vcal_1\to\Vcal_2} =\sup_{0\ne v\in\Vcal_1} \|L v\|_2 / \|v\|_1$.  We can write the right side of \eqref{BFestHilb} in a cluttered but precise manner as ``$\|V\|_{l^2 \to \Xcal} \|V^{-1}\|_{\Xcal \to l^2} \|(B-A)x\|_{\Xcal}$.''  

\begin{proof}[Proof of theorem \ref{thm:BFH}]   If $\mu\in\sigma(A)$ there is nothing to prove.  Otherwise, $\left(A-\mu I\right)^{-1} = V \left(\Lambda-\mu I\right)^{-1} V^{-1}$ so that $\left(\left\|\left(\Lambda-\mu I\right)^{-1}\right\|\right)^{-1} \le \|V\|\, \|V^{-1}\| \, / \, \left\|\left(A-\mu I\right)^{-1}\right\|$.
Note $(\Lambda-\mu I)^{-1}$ is diagonal and bounded with
    $$\|(\Lambda-\mu I)^{-1}\| \le \sup_j |\lam_j-\mu|^{-1} = \max_{\lam\in\sigma(A)} |\lam-\mu|^{-1} = \left(\min_{\lam\in\sigma(A)} |\lam-\mu|\right)^{-1},$$
where $\{\lam_i\}$ are the diagonal entries of $\Lambda$.  Thus
\begin{equation}\label{estwres}
\min_{\lam\in\sigma(A)} |\lam-\mu| \le \frac{\|V\| \|V^{-1}\|}{\left\|\left(A-\mu I\right)^{-1}\right\|}.
\end{equation}

On the other hand, if $\mu$ is an eigenvalue of $B$ and $x\in {H^1}$ satisfies $Bx=\mu x$ and $\|x\|=1$, then $(B-A)x=-(A - \mu I)x$, or $x=-\left(A-\mu I\right)^{-1} (B-A) x$.  Taking norms, it follows that $1\le \|(A-\mu I)^{-1}\| \|(B-A)x\|$.  Combining this with inequality \eqref{estwres} we get estimate \eqref{BFestHilb}.\end{proof}

Recalling lemma \ref{lem:eiginterp} we know that there exist polynomials very close to the eigenfunctions of $\hU$.  This motivates the following corollary to theorem \ref{thm:BFH}.

\begin{cor}\label{cor:BFHwp}  Assume the hypotheses of theorem \ref{thm:BFH}.  Additionally, suppose $p\in\Xcal$ satisfies $\|x-p\|<\eps$.  Then
$$\min_{\lam\in\sigma(A)} |\mu-\lam| \le \|V\| \|V^{-1}\| \Big[\eps (\|B\|+\|A\|) + \|(B-A)p\| \Big].$$
\end{cor}

\begin{proof}  Note $\|(B-A)x\| \le \|(B-A)(x-p)\| + \|(B-A)p\| \le (\|B\|+\|A\|) \|x-p\| + \|(B-A)p\|$. \end{proof}

\subsection{Approximating the eigenvalues of $U$}  \label{subsect:approxeig}  Our strategy for estimating the size of the difference between the computable eigenvalues of $U_N$ (actually, $\hU_N$) and the desired eigenvalues of $U$ (actually, $\hU$) follows the outline: 
\renewcommand{\labelenumi}{(\emph{\roman{enumi}})} 
\begin{enumerate}
\item diagonalize $\hU_N=\hV\hat\Lambda \hV^{-1}$ by numerically diagonalizing $U_N=V \Lambda V^{-1}$;
\item estimate $\|\hV\|,\|\hV^{-1}\|$ by noting $P_N V = \trefoneline{\tT_0}{\tT_N} \Gamma$ (see lemma \ref{lem:diagUm}) and finding the matrix $\Gamma$ numerically; note $\|\hU_N\|_{H^1} \le |\lam_1| \|\hV\| \|\hV^{-1}\|$ if $\lam_1$ is the largest eigenvalue of $U_N$;
\item also using $U_N$, approximately solve consecutive initial value problems with history functions $\tT_0, \tT_1, \dots, \tT_N$ (i.e.~the ${H^1}$-normalized Chebyshev polynomials in subsection \ref{subsect:Hone}) and record \emph{a posteriori} estimates from theorem \ref{thm:apostivp};
\item use lemma \ref{lem:eiginterp} to show that an eigenvector of $\hU$ with eigenvalue away from zero is well-approximated by a polynomial $p$ of degree $k\le N$, and thus by a linear combination of $\tT_0, \dots, \tT_k$;
\item (this is Theorem \ref{thm:main}): use corollary \ref{cor:BFHwp} with $\Xcal={H^1}$, $A=\hU_N$, and $B=\hU$; estimate norm $\|\hU\|$ from lemma \ref{lem:Uhatnorm}; bound $\|(B-A)p\|$ by estimates in (\emph{iii}); conclude that eigenvalues of $U$ are close to the computed eigenvalues of $U_N$.\end{enumerate}

\subsection{An \emph{a posteriori} theorem for the eigenvalues of $\hU$}\label{subsect:aposteig}  First we show in detail how diagonalizing $U_N$ diagonalizes $\hU_N$.

\begin{lem}  \label{lem:diagUm} Suppose $U_N = V \Lambda V^{-1}$ with $\Lambda=(\lam_j)_{j=1}^{N+1}$ diagonal.  Let $v_k$ be the $k$th column of $V$, so that $U_N v_k = \lam_k v_k$, and define $p_k=P_N v_k$.   Expand $p_k$ in ${H^1}$-normalized Chebyshev series $p_k(t) = \sum_{j=1}^{N+1} \Gamma_{jk} \tT_{j-1}(t)$, thus defining an invertible matrix $\Gamma\in\CC^{(N+1)\times (N+1)}$.

Denote by $a=(a_1,a_2,\dots)$ an element of $l^2$.  Let $\hV:\l^2 \to {H^1}$ be the operator 
    $$\hV a \equiv \sum_{k=1}^{N+1} a_k p_k + \sum_{k>N+1} a_k \tT_{k-1} = \sum_{j=1}^{N+1} \left(\sum_{k=1}^{N+1} \Gamma_{jk} a_k\right) \tT_{j-1} + \sum_{k>N+1} a_k \tT_{k-1}.$$
Then $\hV$ is bounded and invertible.  Furthermore, $\|\hV\|\le \left(\|\Gamma\|_2^2 + 1\right)^{1/2}$ and $\|\hV^{-1}\|\le \left(\|\Gamma^{-1}\|_2^2 + 1\right)^{1/2}$, where $\|\cdot\|_2$ is the matrix $2$-norm on $\CC^{(N+1)\times (N+1)}$.  

Finally, from the diagonal matrix $\Lambda$, define $\hat\Lambda \in \Lop(l^2)$ by $(\hat\Lambda a)_j = \lam_j a_j$ if $j\le N+1$ while $(\hat\Lambda a)_j=0$ otherwise.  (Thus $\hat\Lambda$ is a diagonal operator of finite rank.)  Then $\hU_N = \hV \hat\Lambda \hV^{-1}$ and also $\|\hU_N\| \le \left(\max |\lam_j|\right) \|\hV\| \|\hV^{-1}\|$.\end{lem}  

\begin{proof}  First, $\Gamma$ is invertible because $V$ and $P_N$ are invertible and $\{\tT_j(t)\}_{j=0}^N$ is a linearly independent set in ${H^1}$.  Next, the inverse of $\hV$ is given by 
    $$\hV^{-1} \left(\sum_{j=1}^\infty \beta_j \tT_{j-1}\right) = \left(\sum_{k=1}^{N+1} (\Gamma^{-1})_{1k}\beta_k,\dots,\sum_{k=1}^{N+1} (\Gamma^{-1})_{N+1,k}\beta_k,\beta_{N+2},\dots\right),$$
a map $\hV^{-1}:{H^1}\to l^2$.  It is easily checked that $\hV^{-1} \hV \alpha = \alpha$ for any $\alpha\in l^2$.

Next,
\begin{align*}
\|\hV a\|_{H^1}^2 &= \sum_{j=1}^{N+1} \left|\sum_{k=1}^{N+1} \Gamma_{jk} a_k\right|^2 + \sum_{j>N+1} |a_j|^2 \le \|\Gamma\|_2^2 \left(\sum_{k=1}^{N+1} |a_k|^2\right) + \sum_{j>N+1} |a_j|^2 \\
    &\le \|\Gamma\|_2^2 \|a\|_{l^2}^2 + \|a\|_{l^2}^2,
\end{align*}
so $\|\hV\|^2 \le \|\Gamma\|_2^2 + 1$, as claimed.  A similar calculation shows $\|\hV^{-1}\|^2 \le \|\Gamma^{-1}\|_2^2 + 1$.

Finally, we show $\hU_N=\hV\hat\Lambda \hV^{-1}$ by action on basis elements of ${H^1}$.  If $j>N+1$ then $\hU_N\tT_{j-1} = P_N U_N E_N \Pi_N \tT_{j-1} = 0$ while $\hV\hat\Lambda\hV^{-1}\tT_{j-1} = \hV\hat\Lambda \delta_j = 0$.  For $j=1,\dots,N+1$, however,
\begin{align}
\hU_N \tT_{j-1} &= P_N U_N E_N \tT_{j-1} = P_N V \Lambda V^{-1} P_N^{-1} \tT_{j-1} = \sum_{k=1}^{N+1} (\Gamma^{-1})_{kj} P_N V \Lambda V^{-1} P_N^{-1} p_k \label{hUmontT} \\
    &= \sum_{k=1}^{N+1} (\Gamma^{-1})_{kj} \lam_k p_k = \sum_{l=1}^{N+1} \Big(\sum_{k=1}^{N+1} \Gamma_{lk} \lam_k (\Gamma^{-1})_{kj}\Big) \tT_{l-1}. \notag\end{align}
We have used the easily checked facts that $\sum_{k=1}^{N+1} (\Gamma^{-1})_{kj} p_k = \tT_{j-1}$ and $V^{-1} P_N^{-1} p_k = e_k$, where $e_k$ is a standard basis element of $\CC^{N+1}$.  On the other hand, by definition of $\hV,\hat\Lambda,$ and $\hV^{-1}$,
\begin{align*}
    \hV \Lambda \hV^{-1} \tT_{j-1} &= \hV \left(\lam_1(\Gamma^{-1})_{1j}, \dots,\lam_{N+1}(\Gamma^{-1})_{N+1,j}, 0,\dots\right) = \sum_{l=1}^{N+1} \Big(\sum_{k=1}^{N+1} \Gamma_{lk} \lam_k (\Gamma^{-1})_{kj}\Big) \tT_{l-1},\end{align*}
exactly the result of \eqref{hUmontT}.\end{proof}

Lemma \ref{lem:diagUm} defines $\hV$ as the product
\begin{equation}\label{hVprodidea}
\hV=\trefonelinefive{\tT_0}{\dots}{\tT_N}{\tT_{N+1}}{\dots} \left(\begin{array}{c|c} \Gamma & 0 \\ \hline 0 & I\end{array}\right)
\end{equation}
where $\Gamma$ is an $(N+1)\times (N+1)$ matrix satisfying
    $$P_N V = \trefoneline{p_1}{p_{N+1}} = \trefoneline{\tT_0}{\tT_N} \Gamma.$$
Here ``$\trefoneline{\cdot}{\cdot}$'' denotes a ``matrix'' whose columns are elements of a function space \cite{BattlesTrefethen}.  The right hand factor in equation \eqref{hVprodidea} is a bounded operator on $l^2$ with ``$I$'' denoting the identity operator on the space spanned by $\{\tT_j\}_{j=N+1}^\infty$.  In any case, if $C$ is the $(N+1)\times (N+1)$ matrix described in equation \eqref{discreteCheb} in section \ref{sect:interp} then
    $$\Gamma = \sqrt{\pi/2}\, \operatorname{diag}\left(\,2^{-1/2},\,2,\,3,\,\dots,\,N+1\right) C V.$$

We can now state the main theorem on numerical approximation of $U$.  It says that the difference between the eigenvalues of the (approximate) monodromy matrix $U_N$ and the exact eigenvalues of the monodromy operator $U$ can be bounded using the machinery developed in this and the last section and also \emph{a posteriori} information for specific initial value problems provided by theorem \ref{thm:apostivp}.

\begin{thm}\label{thm:main}  Let $N\ge 1$.  Suppose $a,b$ in DDE \eqref{DDEscalar} are analytic on $I=[-1,1]$ with regularity ellipse $E\supset I$.  Let $e^\eta = S+s>1$ be the sum of semiaxes of $E$.  Let $A_E=\max_{z\in E} \left|\int_{-1}^z a\right|$, $B_E=\max_{z\in E} \left|\int_{-1}^z b\right|$.

Suppose $\mu\in\CC$ is an eigenvalue of $U$ (and of $\hU$) such that $|\mu|\ge\delta>0$.  Let
    $$\eps_k = \left(8 (\sinh\eta)^{-1}\, e^{A_E+ B_E/\delta}\right) \, k \, e^{-k\eta}$$
for integers $k\ge 1$.

Assume $U_N = V \Lambda V^{-1}$ is diagonalizable, with $\Lambda=(\lam_j)$ a diagonal matrix and eigenvalues $\{\lam_j\}_{j=1}^{N+1}$ ordered by decreasing magnitude.  This gives a diagonalization $\hU_N = \hV \hat \Lambda \hV^{-1}$ as in lemma \ref{lem:diagUm}.  Let $\cond(\hV)=\|\hV\| \|\hV^{-1}\|$.

Suppose $\|\hU\tT_j - \hU_N \tT_j\|_{H^1} \le \nu_j$ for $j=0,\dots,N$ and let $\xi_k=\Big(\sum_{j=0}^k \nu_j^2\Big)^{1/2}$ for $k=1,\dots,N$.  Define
\begin{equation}\label{omegadefn}
\omega_k = \eps_k \left(\|\hU\| + |\lam_1| \cond(\hV)\right) + (1+\eps_k) \xi_k.
\end{equation}
Then 
\begin{equation}\label{mainest}
\min_{j=1,\dots,N+1} |\mu-\lam_j| \le \min \left\{\omega_1,\dots,\omega_N\right\} \cond(\hV).
\end{equation}
\end{thm}

\smallskip
\begin{rem}  Note that one computes $\nu_j$ by solving the ODE initial value problems
    $$\dot y = a(t) y + b(t) \tilde T_j(t), \qquad y(-1) = \tilde T_j(1),$$
and then by using theorem \ref{thm:apostivp}, corollary \ref{cor:apost}, and lemma \ref{lem:ConethenHcal}.  Also, $\|\hU\|$ is estimated in lemma \ref{lem:Uhatnorm} in terms of $\|a\|_\infty$, $\|b\|_\infty$, and $C_a$; the last constant bounds the fundamental solution $\Phi_a$. \end{rem}

\smallskip
\begin{proof}  Let $x\in{H^1}$ be a normalized eigenvector associated to $\mu$: $\hU x = \mu x$, $\|x\|_{H^1}=1$.  By lemma \ref{lem:eiginterp}, for each $k=1,\dots,N$ there is $q_k\in\Po_k$ such that $\|x-q_k\|_{H^1}\le \eps_k$.

Apply corollary \ref{cor:BFHwp} with $\Xcal={H^1}$, $A=\hU_N$, $B=\hU$, and $p=q_k$.  Also, use lemma \ref{lem:diagUm} to diagonalize $\hU_N=\hV \hat\Lambda \hV^{-1}$, with $\sigma(\hU_N)=\{\lam_j\}_{j=1}^{N+1}\cup\{0\}$ and find
\begin{equation}\label{estfirstform}
\min_{j=1,\dots,N+1} |\mu-\lam_j| \le \|\hV\| \|\hV^{-1}\| \left(\eps_k (\|\hU\|+\|\hU_N\|) +\|(\hU-\hU_N)q_k\|_{H^1}\right).
\end{equation}

On the other hand, $q_k = \sum_{j=0}^k \alpha_{kj} \tilde T_j(t)$ so
\begin{align*}
\|(\hU-\hU_N) q_k\|_{H^1} &\le \sum_{j=0}^k |\alpha_{kj}| \|(\hU-\hU_N)\tilde T_j\|_{H^1} \le \sum_{j=0}^k |\alpha_{kj}| \,\nu_j \\
    &\le \|q_k\|_{H^1} \, \xi_k  \le (\|x\|_{H^1} + \|q_k-x\|_{H^1}) \xi_k \le (1+\eps_k) \xi_k\end{align*}
by Cauchy-Schwarz.  From \eqref{estfirstform} and lemma \ref{lem:diagUm} we see that
    $$\min_{j=1,\dots,N+1} |\mu-\lam_j| \le \cond(\hV) \left(\eps_k \left(\|\hU\| + |\lam_1|\cond(\hV)\right) + (1+\eps_k) \xi_k\right).$$
as desired.\end{proof}

\subsection{Using the theorem in an example}\label{subsect:usethm}  The theory of the last section is complicated.  It remains to show by example that the bounds given there are practical.  

It is worthwhile to give heuristics about the behavior of the constants $\eps_k$, $\nu_j$, $\xi_k$, and $\omega_k$.  First, for fixed analytic functions $a(t)$ and $b(t)$, $\eps_k$ is monotonic decreasing and decays exponentially with rate arbitrarily close to $\frac{1}{S+s} = e^{-\eta} < 1$.  Second, $\nu_j$ is (roughly) increasing as $j$ increases.  Examples show it is commonly within a few orders of magnitude of machine precision $\eps_m$ for $j$ up to approximately $0.7N$, supposing $N$ is sufficiently large to achieve $\eps_m$ accuracy for interpolation of the coefficient functions $a(t),b(t)$.  The meaning of ``few orders'' is dependent on stiffness, however---see subsection \ref{subsect:ivpexamples}.  Third, $\xi_k$  accumulates the $\nu_j$ so $\xi_k$ is also a few orders of magnitude above $\eps_m$.  As $N$ is typically in the range $10<N<300$, accumulation is small and $\xi_k$ is at most two orders of magnitude above the maximum $\nu_j$ for $j\le k$.  On the other hand, $\nu_j$ will be large for $j\approx N$ because $\hU_N \tilde T_j$ is a poor approximation of $\hU \tilde T_j$ in that case, and thus $\xi_k$ can become $O(1)$ for $k\approx N$.

We conclude that $\omega_k$ typically decreases as $\eps_k$ decreases until $\eps_k \left(\|\hU\|+|\lam_1|\cond(\hV)\right) \approx \xi_k$ at which point the minimum value $\omega$ is achieved.

Note, also, that $\cond(\hV)$ is generally an increasing function of $N$ which seems to be of size $10^3$ to $10^6$ for $40\le N \le 300$ as we see in the following example.  In fact, figure \ref{fig:constantsdecay}\textbf{b} suggests that $\cond(\hV)$ grows subexponentially in the example below.  Such behavior is evidently desirable for our application.

\begin{example} Consider the DDE,
\begin{equation}\label{exDDE}
\dot y = a y + \left(b + \sin(3 \pi t)\right) y(t-2)
\end{equation}
which is equation \eqref{introex} of the introduction (section \ref{sect:intro}).  Fix parameter values $(a,b)=(-1.1,1)$.  We find numerically that $\sigma(U_N)=|\lam_1|=0.9369$ for the monodromy matrix, for all $N$ greater than about $15$.  We want to confirm this approximation and thereby confirm stability for this parameter pair.

Though $b(t)=b + \sin(3 \pi t)$ is entire, we nonetheless must choose a regularity ellipse $E$ to find the constants in theorem \ref{thm:main}.  In particular, the oscillation of $b(t)$ for $t\in I$ is related to the growth of $|b(z)|$ for $z$ with increasing imaginary part.  In fact, note that
    $$B_E = \max_{z\in E} \left|\int_{-1}^z b + \sin(3\pi\zeta)\,d\zeta\right| \le |b|(1+S) + (3\pi)^{-1} (\cosh(3\pi s)+1).$$
Thus one wants to choose the smaller semiaxis $s$ to be relatively small.  It turns out that the choices $s=0.5$ and $S=\sqrt{1+s^2}\approx 1.1$ are reasonable.  

Figure \ref{fig:constantsdecay}\textbf{a} shows the growth of $\xi_k$, the decay of $\eps_k$, and the (non-monotonic) behavior of $\omega_k$ when $N=184$.  Note that $\omega=\min \omega_k$ occurs at $k\approx 160$.  The constants $\eps_k$ in the statement of theorem \ref{thm:main} are quite large when $k$ is small.  In addition to being sensitive to the regularity ellipse for $b(t)$, as just described, the lower bound $\delta$ on interesting eigenvalues controls the size of $N$ necessary to achieve a given accuracy.  Here we let $\delta=0.2$, that is, we consider the eigenvalues $\mu$ of $U$ such that $|\mu|\ge 0.2$.

\begin{figure}[ht]
\regfigure{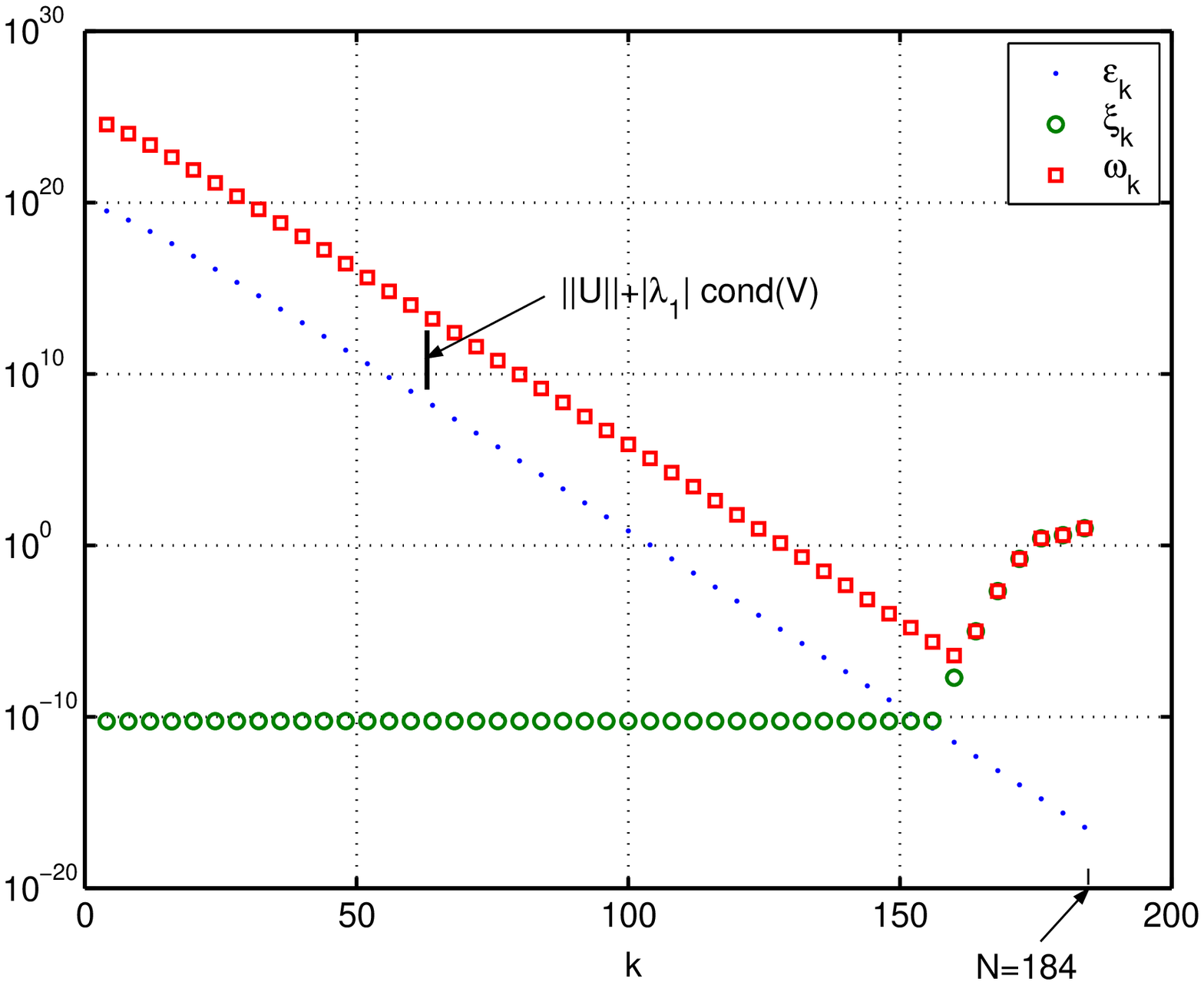}{2.4}\quad\regfigure{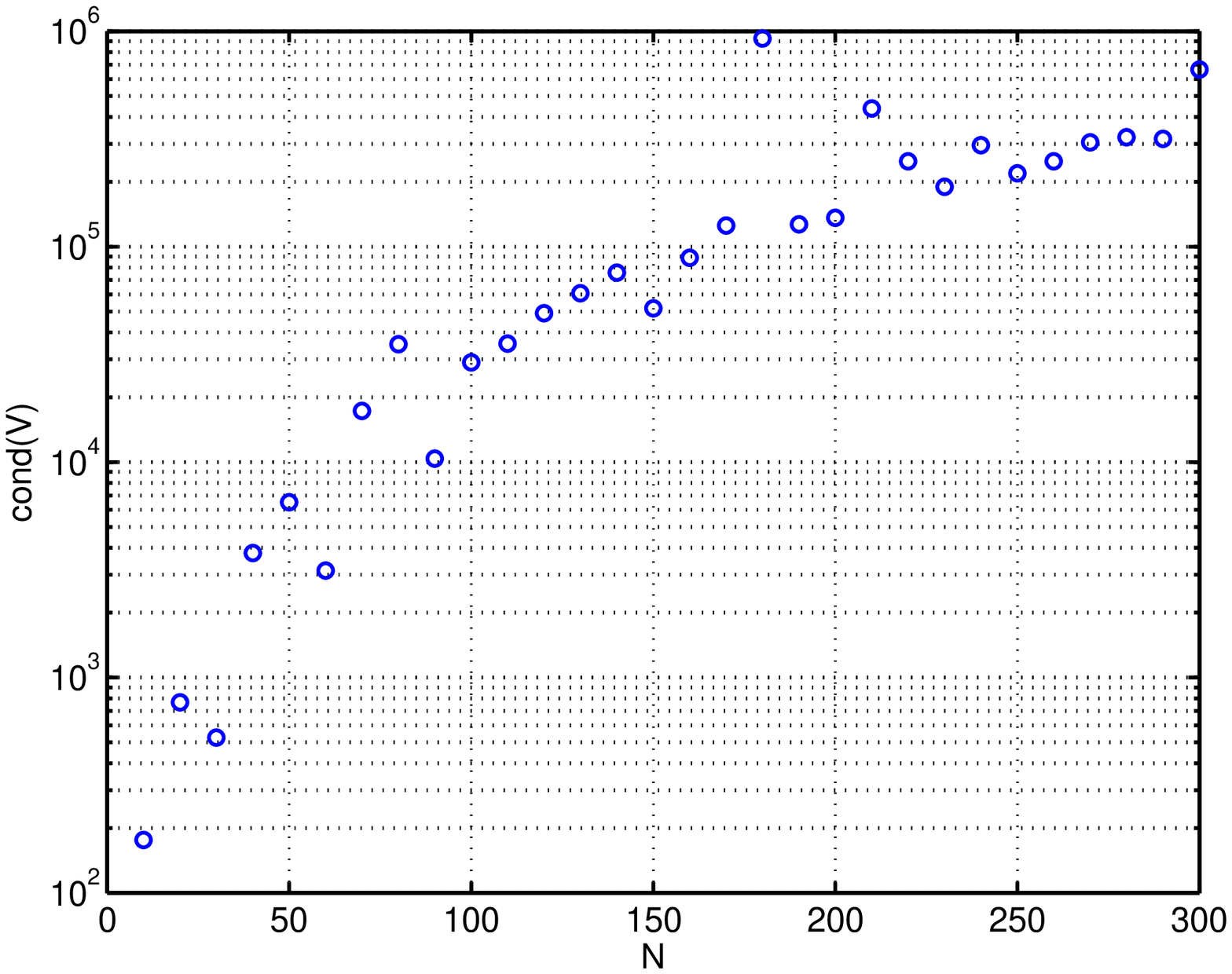}{2.4}
\caption{\textbf{(a)}  Behavior of constants $\eps_k$, $\xi_k$, and $\omega_k$ (see text) in theorem \ref{thm:main}, for DDE \eqref{exDDE} and $N=184$.  \textbf{(b)} $\cond(\hV)$ as a function of $N$.} \label{fig:constantsdecay}
\end{figure}

The result of applying theorem \ref{thm:main} is that the error radius is approximately $0.0434$ so we have proven that $|\mu|<1$ for all eigenvalues of $U$ as shown in figure \ref{fig:introeigs} of section \ref{sect:intro}.  As $N$ increases above $N=184$, furthermore, we see that the error radius decreases exponentially (figure \ref{fig:eradexpdecay} in section \ref{sect:intro}).
\end{example}

\section{Generalization to systems}  \label{sect:systems}

\subsection{Bounds on fundamental solutions}  Although there are other technical complications, generalization of the previous results to systems is mainly distinguished by diminished understanding of the fundamental solution to the homogeneous ODE problem.

Recall, in particular, that theorem \ref{thm:apostivp}, corollary \ref{cor:apost} and lemma \ref{lem:Uhatnorm} for the scalar case used the bound $|\Phi_a(t)|\le C_a = \exp\left(\int_{-1}^1 \max\{\Re a(s),0\}\,ds\right)$ on the fundamental solution $\Phi_a(t)$ to $\dot y(t) = a(t) y(t)$.  On the other hand, lemma \ref{lem:eiginterp} required a uniform bound on the analytical continuation $\Phi_\lam(z), z\in E$, of $\Phi_\lam(t)$ solving $\dot x(t) = (a(t)+b(t)/\lam) x(t)$.

We start with the existence and analytic continuation of the fundamental solution.

\begin{lem}\label{lem:fundanalbound}  Suppose $A(z)\in \CC^{d\times d}$ is (entry-wise) analytic on a convex open set $E\supset I=[-1,1]$.  Then there is a unique function $\Phi_A(z)$, analytic on $E$, satisfying
\begin{equation}\label{zintODE}
\Phi_A(z) = I_d + \int_{-1}^z A(\zeta) \Phi_A(\zeta)\,d\zeta.
\end{equation}
If $|A(z)|\le \alpha$ for $z\in E$ then $|\Phi_A(z)|\le e^{\alpha |z+1|}$ for $z\in A$.  Furthermore $\dot\Phi_A(t) = A(t) \Phi_A(t)$ for $t\in I$. \end{lem}

\begin{proof}  For $z\in E$ choose a path $\gamma(s)=-1+(z+1) s/r$ for $r=|z+1|$ and $0\le s \le r$.  Picard iteration on \eqref{zintODE}, for $z\in E$, is
    $$\Phi_q(z)=I_d + \frac{z+1}{r} \int_{0}^r A(\gamma(s)) \Phi_{q-1}(\gamma(s))\,ds, \quad \Phi_0(z)=I_d.$$
If $\eta_q(z)=|\Phi_q(z)-\Phi_{q-1}(z)|$ then induction gives $\eta_q(z) \le \alpha^q r^q/q!$, so $\{\Phi_q(z)\}$ is a Cauchy sequence of analytic functions on $E$.  Thus $\Phi_A(z)= \lim_q \Phi_q(z)$ is analytic.  Equation \eqref{zintODE} implies $\dot \Phi_A(t)=A(t) \Phi_A(t)$ for $t\in I$.  Finally,
\begin{align*}
|\Phi_A(z)| &= \lim_q \left|I_d + \sum_{j=1}^d \Phi_j(z)-\Phi_{j-1}(z)\right| \le \sum_{j=0}^\infty \alpha^j r^j/j! = e^{\alpha r}. \end{align*}
\end{proof}

One can easily show that, in addition,

\begin{cor}  Fix $s\in I$.  Let $\Omega(t)=\Phi_A(t) \Phi_A(s)^{-1}$ for $s\le t\le 1$.  Then $\Omega$ satisfies $\dot\Omega(t) = A(t)\Omega(t)$, $\Omega(s)=I_d$ and furthermore $|\Omega(s)|\le e^{\tilde\alpha(t-s)}$ if $|A(\tau)|\le \tilde\alpha$ for $\tau\in[s,t]$.
\end{cor}

There is, however, a more refined source of bounds on fundamental solutions.  Namely, one can use the collocation algorithm (below) to approximate the fundamental solution and then use \emph{a posteriori} estimates from theorem \ref{thm:apostivp} to bound the fundamental solution.  There is a ``bootstrapping'' aspect to this: one must have some bound on the fundamental solution in order to compute the \emph{a posteriori} estimates which leads to an improved bound.  This is the content of lemma \ref{lem:collocboundfund} below.

\subsection{The collocation algorithm for initial value problems}\label{subsect:collocalgsys}  We now describe the collocation algorithm for initial value problems of the form
\begin{equation}\label{odesys}
\dot\by(t)=A(t)\by(t)+\bu(t), \qquad \by(-1)=\by_0,
\end{equation}
for $\by(t), \bu(t), \by_0 \in\CC^d$ and $A(t)$ a $d\times d$ matrix.

As before, fix the interval $I=[-1,1]$ for $t$.  Define the space of $\CC^d$-valued continuous functions $\bC(I)=C(I)\otimes \CC^d$ so that $\bbf=(f_1,\dots,f_d)^\top$ is a $d\times 1$ column vector of functions $f_j\in C(I)$.  For $\bbf\in\bC(I)$ let $\|\bbf\|_\infty = \max_{t\in I} |\bbf(t)|$ where ``$|\cdot|$'' is the Euclidean norm on $\CC^d$, that is, $|\bv|^2 = \sum_{j=1}^d |v_j|^2$.  Note $|\cdot|$ induces a norm on $d\times d$ matrices; we also denote this norm ``$|\cdot|$.''  Continuous matrix-valued functions $A(t)=\left(a_{ij}(t)\right)$ by definition satisfy $a_{ij}\in C(I)$ for all $i,j$.  For such $A$ define the norm $\|A\|_\infty = \max_{t\in I} |A(t)|$.

Let $\bPo_N$ be the space of $\CC^d$-valued polynomials of degree at most $N$.  

On the function spaces $\bC(I)$ and $\bPo_N$ we have collocation operators.  In particular, evaluation at $N+1$ Chebyshev collocation points $\Co_N=\{t_0,\dots,t_N\} = \{\cos(j\pi/N)\}$ to produces a vector in $\CC^l$ where $l=d(N+1)$:
\begin{equation}\label{collocevalonsys}
\bEcal_N \bbf = \left(f_1(t_0),\dots,f_d(t_0),f_1(t_1),\dots,f_d(t_1), \dots,\dots,\dots, f_1(t_N),\dots,f_d(t_N)\right)^\top.
\end{equation}
(From now on, we order the scalar components of $\CC^l$ consistently with the output of $\bEcal_N$.)  Restricting $\bEcal_N$ to polynomials gives $\bE_N:\bPo_N\to\CC^l$.  The inverse of $\bE_N$, namely $\bP_N:\CC^l\to\bPo_N$, creates a $\CC^d$-valued polynomial from collocation values.  The interpolation operator is $\bI_N = \bP_N \bEcal_N : \bC(I) \to \bPo_N$.  That is, if $\bp=\bI_N \bbf$ for $\bbf\in\bC(I)$ then $\bp$ is a $\CC^d$-valued polynomial of degree $N$ such that $\bp(t_j)=\bbf(t_j)$.

The ``Chebyshev (spectral) differentiation matrix'' is the map
    $$\bD_N = \bE_N \circ \frac{d}{dt} \circ \bP_N : \CC^l \to \CC^l.$$
If $D_N$ is the scalar differentiation matrix (section \ref{sect:intro}) then $\bD_N =D_N\otimes I_d$.  In \Matlab, using \mtt{cheb.m} from \cite{Trefethen}, if \mtt{D=cheb(N)} then $\bD_N =$ \mtt{kron(D,eye(d))}.

Recall that in subsection \ref{subsect:collocalg} we defined matrices $\hD_N$ and $\hM_a$ and also a vector $\hu$ in order to write the collocation algorithm as a matrix calculation.  Here we define
\begin{gather*}
\hbD_N = \hD_N\otimes I_d, \quad \hM_A = \begin{pmatrix}A(t_0) & & & \\ & \ddots & & \\ & & A(t_{N-1}) & \\ & & & 0_d \end{pmatrix}, \quad \text{and} \quad \hbu=\begin{pmatrix} \bu(t_0) \\ \vdots \\ \bu(t_{N-1}) \\ \by_0\end{pmatrix},
\end{gather*}
where $I_d,0_d$ are the $d\times d$ identity and zero matrices, respectively.  Thus $\hbD_N,\hM_A$ are $l\times l$ matrices and $\hbu\in\CC^l$.

\begin{lem}[the collocation algorithm for systems] \label{lem:collocalgsys} Let $N\ge 1$ and $t_j=\cos(\pi j/N)$, $j=0,\dots,N$.  For $A(t)$ a continuous $d\times d$ matrix-valued function of $t\in I$, $\bu(t)\in \bC(I)$, and $\by_0\in\CC^d$, the following are equivalent:\begin{itemize}
\item $\bp\in\bPo_N$ satisfies
    $$\dot\bp(t_j) = A(t_j) \bp(t_j) + \bu(t_j), \quad 0\le j \le N-1, \qquad \text{and } \quad \bp(-1)=\by_0;$$
\item$\bv\in\CC^l=\CC^{d(N+1)}$ satisfies
    $$\hbD_N \bv = \hM_A \bv + \hbu.$$\end{itemize}
The equivalence is $\bp=\bP_N \bv$ and $\bv=\bE_N \bp$.\end{lem}

\begin{proof} Immediate from the definitions.\end{proof}

For the collocation algorithm we have the following \emph{a posteriori} estimates:

\begin{thm}\label{thm:apostsys}  Let $N\ge 1$.  Suppose $A(t)$ is a continuous $d\times d$ matrix-valued function of $t\in I$, $\bu(t)\in \bC(I)$, and $\by_0\in\CC^d$.  Suppose
    $$\left|\Phi_A(t) \Phi_A(s)^{-1}\right| \le C_A \text{ for all } -1\le s\le t \le 1.$$
Suppose $\by\in \bC(I)$ satisfies the initial value problem \eqref{odesys}.  Let $\bp\in \bPo_N$ be the $\CC^d$-valued polynomial described by lemma \ref{lem:collocalgsys} and let $R_{\bp} = \dot\bp(-1) - A(-1) \by_0 - \bu(-1)$.  Then
\begin{equation}\label{apostsys}
\|\by-\bp\|_\infty \le 2 C_A \Big[\|A \bp - \bI_N(A\bp)\|_\infty + \|\bu-I_N(\bu)\|_\infty + |R_{\bp}|\Big]
\end{equation}
and
\begin{equation}\label{apostsysdot}
\|\dot\by-\dot\bp\|_\infty \le \left(2 \|A\|_\infty C_A + 1\right)\, \Big[\|A \bp - \bI_N(A\bp)\|_\infty + \|\bu-I_N(\bu)\|_\infty + |R_{\bp}|\Big].
\end{equation}\end{thm}

\begin{proof}  Let $\bq = I_N(A\bp)$ and $\bw=I_N(\bu)$.  Since $\br = \dot \bp-\bq-\bw \in \bPo_N$ and $\br(t_j)=0$ for $j=0,\dots,N-1$, it follows that there is $\bv\in\CC^d$ such that $\br=\bv l_N$.  (Recall $l_N(t)$ is a polynomial defined in lemma \ref{lem:lNtotrig}.)  Evaluating at $t=-1$, we find $R_{\bp} = \bv (-1)^N N 2^{2-N}$ so $|\bv|=|R_{\bp}|\, 2^{N-2}/N$.

On the other hand,
\begin{equation}\label{diffofsys}
\dot\by - \dot\bp = A(\by-\bp) + (A\bp-\bq) + (\bu-\bw) - \bv l_N, \quad (\by-\bp)(-1) = 0,
\end{equation}
so
    $$\by(t)-\bp(t) = \int_{-1}^t \Phi_A(t) \Phi_A(s)^{-1} \Big[(A\bp(s)-\bq(s)) + (\bu(s)-\bw(s)) - \bv l_N(s)\Big]\,ds.$$
Taking $|\cdot|$ norms and using lemma \ref{lem:lNtotrig} to find $\|\bv l_N\|_\infty \le |R_{\bp}|$,
    $$|\by(t)-\bp(t)| \le 2 C_A \Big[\|A \bp - \bq\|_\infty + \|\bu-\bw\|_\infty + |R_{\bp}|\Big],$$
giving \eqref{apostsys}.  For estimate \eqref{apostsysdot}, using equation \eqref{diffofsys} and lemma \ref{lem:lNtotrig}, we find that \eqref{apostsys} implies \eqref{apostsysdot}.\end{proof}

We now return to the computation of bounds on fundamental solutions.  For the initial value problem \eqref{odesys}, first calculate a typically poor bound $\left|\Phi_A(t)\Phi_A(s)^{-1}\right| \le e^{2\|A\|_\infty} =C_A$.  Then approximate both the fundamental solution $\Phi_A(t)$ and its inverse $\Phi_A(t)^{-1}$ by the collocation algorithm.  Then use the \emph{a posteriori} bounds on the error in computing the fundamental solution to bound the norm of the approximate fundamental solution.  Repeat as necessary.

The following lemma gives a precise outline.  Recall that $\Psi_A(t)=\Phi_A(t)^{-\top}$ satisfies the ``adjoint equation'' $\dot\Psi_A(t) = - A(t)^\top \Psi_A(t)$, $\Psi_A(-1)=I_d$.

\begin{lem}\label{lem:collocboundfund}  Consider the following initial value problems:
\begin{equation}\label{fundivps}
\dot \by_s(t)=A(t) \by_s(t), \qquad \by_s(-1)=\be_s, \qquad s=1,\dots,d,
\end{equation}
where $\{\be_s\}$ is the standard basis for $\CC^d$.  Suppose that for $s=1,\dots,d$, $\|\by_s-\bp_s\|_\infty \le \nu_s$ where $\bp_s$ is the collocation algorithm approximation of $\by_s$.  

Note $\Phi_A(t) = \trefoneline{\by_1(t)}{\by_d(t)}$ is the fundamental solution to $\dot \by = A(t) \by$.  Let $\Phi_N(t) = \trefoneline{\bp_1(t)}{\bp_d(t)}$.  If $\xi^2=\sum_{s=1}^d \nu_s^2$ then $|\Phi_A(t)-\Phi_N(t)| \le \xi$ for all $t\in I$.  If, furthermore, $\Psi_A(t)$ is the fundamental solution to the adjoint equation $\dot \bz=-A(t)^\top \bz$, and if $\Psi_N(t)$ is the collocation approximation of $\Psi_A(t)$, and if $\left|\Psi_A(t)-\Psi_N(t)\right|\le \omega$ for $t\in I$, then
    $$\left|\Phi_A(t)\Phi_A(s)^{-1}\right| \le \left(\xi + \|\Phi_N\|_\infty\right) \left(\omega + \|\Psi_N\|_\infty\right).$$
\end{lem}

\begin{proof}  
\begin{align*}
|\Phi_A(t)-\Phi_N(t)| &= \max_{\bu\in\CC^d,\,|\bu|=1} \left|\left(\Phi_A(t)- \Phi_N(t)\right)\bu\right| = \max_{|\bu|=1} \left|\sum_{s=1}^d u_s \left(\by_s(t)-\bp_s(t)\right)\right| \\
    &\le \max_{|\bu|=1} \left(\sum_{s=1}^d |u_s|^2\right)^{1/2} \left(\sum_{s=1}^d |\by_s(t)-\bp_s(t)|^2\right)^{1/2} \le \xi
\end{align*}
by the Cauchy-Schwarz inequality.  To conclude, note
\begin{align*}
|\Phi_A(t)\Phi_A(s)^{-1}| &\le |\Phi_A(t)|\, |\Psi_A(s)| \\
        &\le \left(\|\Phi_A - \Phi_N\|_\infty + \|\Phi_N\|_\infty\right) \left(\|\Psi_A-\Psi_N\|_\infty + \|\Psi_N\|_\infty\right).
\end{align*}\end{proof}

\begin{example}\label{example:getboundMathieu}  The second order ODE $\ddot x + \dot x + (10+9\cos(\pi t))x=0$, a relatively stiff damped Mathieu equation, corresponds to $A(t)=\begin{pmatrix} 0 & 1 \\ -10-9\cos(\pi t) & -1\end{pmatrix}$.  It follows that $|A(t)|,|A(t)^\top| \le |A(t)|_F =\sqrt{2+(10+9\cos(\pi t))^2} \le \sqrt{363}$ for all $t$, where $|\cdot|_F$ denotes the Frobenius norm \cite{BurdenFaires}.  Thus $C_1 = e^{2\|A\|_\infty} \approx 3.5387 \times 10^{16}$ is an \emph{a priori} bound on $|\Phi_A(t)\Phi_A(s)^{-1}|$.  We use the collocation algorithm with $N=50$ to find $\Phi_N(t),\Psi_N(t)$ approximating $\Phi_A(t),\Psi_A(t)$.  The \emph{a posteriori} estimates from theorem \ref{thm:apostsys} are computed using $C_A=C_1$ and, as in lemma \ref{lem:collocboundfund}, we find $|\Phi_A(t)\Phi_A(s)^{-1}| \le C_2 = 2.7117 \times 10^9$.  This is a significant improvement, but also we can now iterate, using $C_A=C_2$ in the \emph{a posteriori} estimates to generate $C_3$, and so on.  The result is a sequence of bounds
    $$|\Phi_A(t)\Phi_A(s)^{-1}| \le 3.5387 \times 10^{16}, 2.7117 \times 10^9, 19.627, 19.587, 19.587, \dots$$
Each number on the right \emph{is} a bound, with an \emph{a priori} argument for the first and \emph{a posteriori} arguments for the remainder.  Evidently they converge superlinearly to about $19.587$.  By looking at (numerical) time-dependent components of the fundamental solutions we see that this is near optimal for the Frobenius norm.  In any case, this improvement by 15 orders of magnitude, and comparable improvements for other examples, makes further error estimation practical, as we will see in the following subsections.\end{example}

\subsection{The monodromy operator for a system of DDEs}  \label{subsect:DDEsystemmono}  Let $\bL^2=L^2\otimes \CC^d$ be the Hilbert space of measurable $\bbf=(f_1,\dots,f_d)^\top : I \to \CC^d$ such that $f_s\in L^2$ for $s=1,\dots,d$.  Recall that $L^2=L^2_T(I)$ is the weighted space defined in subsection \ref{subsect:Hone}: $\|f\|_{L^2}^2 = \int_{-1}^1 |f(t)|^2 (1-t^2)^{-1/2}\,dt$.  Similarly, define $\bH^1 = H^1 \otimes \CC^d$, a Sobolev space in which first derivatives are controlled.  Note that $\|\bbf\|_{\bH_1}^2 = \sum_{s=1}^d \|f_s\|_{H^1}^2$ and that lemma \ref{lem:pointwiseeval} says $|\bbf(t)|\le 0.9062\|\bbf\|_{\bH^1}$ for fixed $t\in I$; our choice of Euclidean norm on $\CC^d$ is important for these latter statements.  Also, from lemma \ref{lem:ConethenHcal} it follows that $\|\bbf\|_{\bH^1}^2 \le 2\pi\left(\|\bbf\|_\infty^2 + \|\dot\bbf\|_\infty^2\right)$.

We address the linear DDE
\begin{equation}\label{DDEsys}
\dot\by(t) = A(t) \by(t) + B(t) \by(t-2)
\end{equation}
for $A,B$ continuous and periodic (with period $T=2$) $d\times d$ matrix-valued functions of $t\in I$.  The monodromy operator for this DDE is $\hbU\in\Lop(\bH^1)$ defined by
    $$(\hbU \bbf)(t) = \Phi_A(t) \left[\bbf(1) + \int_{-1}^t \Phi_A(s)^{-1} B(s) \bbf(s)\,ds\right]$$
and $\by_n = \hbU^n \bbf$ solves by method-of-steps the initial value problem consisting of \eqref{DDEsys} and $\by(t)=\bbf(t+2)$, $t\in I$, $\bbf\in\bH^1$.

If $\hbU \bx =\lam \bx$ for $\lam\in\CC$ then $\bz=\lam\bx$ solves $\dot\bz = A\bz + B \bx$.  Thus if $\lam\ne 0$ then $\dot\bx = \left(A + \lam^{-1} B\right)\bx$.  If $\Phi_\lam(t)$ is the solution to $\dot\Phi_\lam = (A+\lam^{-1} B)\Phi_\lam$ and $\Phi_\lam(-1)=I_d$ then $\bx(t) = \Phi_\lam(t) \bx(-1)$.

The next lemmas are systems analogs of lemmas \ref{lem:Uhatnorm} and \ref{lem:eiginterp}.

\begin{lem}\label{lem:Uhatnormsys}  Suppose $\left|\Phi_A(t) \Phi_A(s)^{-1}\right| \le C_A$ for all $-1\le s\le t \le 1$.  Let $a^2=1+\|A\|_\infty^2$ and $c=0.9062$.  Then
\begin{equation}\label{hbUnormest}
\|\hbU\|_{\bH^1} \le \sqrt{2\pi d} \left(c a C_A + \|B\|_\infty \left(c^2 +\pi a^2 C_A^2/2\right)^{1/2}\right).
\end{equation}
\end{lem}

\begin{proof}  Suppose $\bbf\in\bH^1$ and let $\bg(t)=\int_{-1}^t \Phi_A(s)^{-1} B(s)\bbf(s)\,ds$.  Note $\|\hbU \bbf\|_{\bH^1} \le \|\Phi_A(t) \bbf(-1)\|_{\bH^1} + \|\Phi_A(t) \bg(t)\|_{\bH^1}$, with the first term bounded using lemma \ref{lem:ConethenHcal}:
\begin{align*}
&\|\Phi_A(t)\bbf(-1)\|_{\bH^1}^2 =\sum_{k=1}^d \|(\Phi_A(t)\bbf(-1))_k\|_{H^1}^2 \\
    &\qquad \le 2\pi \sum_{k=1}^d \|(\Phi_A(t)\bbf(-1))_k\|_\infty^2 + \|(A(t) \Phi_A(t)\bbf(-1))_k\|_\infty^2 \le 2\pi d \,C_A^2 a^2 |\bbf(-1)|^2.
\end{align*}

On the other hand, $\frac{d}{dt} \left(\Phi_A(t)\bg(t)\right) = A(t) \Phi_A(t) \bg(t) + B(t)\bbf(t)$ so by lemma \ref{lem:ConethenHcal},
    $$\|\Phi_A(t)\bg(t)\|_{\bH^1}^2 = 2\pi \sum_{k=1}^d \|(\Phi_A(t)\bg(t))_k\|_\infty^2 + \left(\|(A(t)\Phi_A(t)\bg(t))_k\|_\infty + \|(B(t)\bbf(t))_k\|_\infty\right)^2.$$
But
\begin{align*}
|(\Phi_A(t)\bg(t))_k| &\le |\Phi_A(t)\bg(t)| \le \int_{-1}^1 \max_{-1\le s\le t\le 1} \left|\Phi_A(t)\Phi_A(s)^{-1}\right|\,|B(s)|\,|\bbf(s)|\,ds \\
    &\le C_A \|B\|_\infty \int_{-1}^1 |\bbf(s)|\,ds \le \sqrt{\frac{\pi}{2}} C_A \|B\|_\infty \|\bbf\|_{\bL^2} \le \sqrt{\frac{\pi}{2}} C_A \|B\|_\infty \|\bbf\|_{\bH^1},
\end{align*}
by the Cauchy-Schwarz inequality with weight $(1-s^2)^{-1/2}\,ds$.  Similarly,
    $$\left|(A(t)\Phi_A(t)\bg(t))_k\right| \le \sqrt{\frac{\pi}{2}} \|A\|_\infty C_A \|B\|_\infty \|\bbf\|_{\bH^1}.$$
On the other hand, $\left|(B(t) \bbf(t))_k\right| \le \|B\|_\infty |\bbf(t)| \le c \|B\|_\infty \|\bbf\|_{\bH^1}$.  Thus
    $$\|\Phi_A(t)\bg(t)\|_{\bH^1}^2 \le 2\pi d \|B\|_\infty^2 \left(\frac{\pi}{2} C_A^2 + \frac{\pi}{2} \|A\|_\infty^2 C_A^2 + c^2\right)\|\bbf\|_{\bH^1}^2,$$
and estimate \eqref{hbUnormest} follows.\end{proof}

\begin{lem}\label{lem:eigsys}  Suppose $A,B$ are analytic $d\times d$ matrix-valued functions of $t\in I$ with common regularity ellipse $E\subset \CC$ which has foci $\pm 1$ and sum of semiaxes $e^\eta = S+s >1$.  Suppose $\hbU\bx=\lam \bx$ for $\lam\ne 0$ and $\|\bx\|_{\bH^1} = 1$.  Let $\Phi_\lam(z)$ be the unique analytic continuation of $\Phi_\lam(t)$ for $z\in E$ and suppose $C_\lam$ is a bound for $|\Phi_\lam(z)|$ if $z\in E$.  If $\bp = \bI_k \bx$ then
    $$\|\bx-\bp\|_\infty \le 8\sqrt{d}\, C_\lam (\sinh\eta)^{-1} k\, e^{-k\eta}.$$
\end{lem}

\begin{proof}  First, $\bx(z)=\Phi_\lam(z) \bx(-1)$ is the analytic continuation of $\bx(t)$, $t\in I$, to $z\in E$.  Furthermore, $|\bx(z)|\le C_\lam |\bx(-1)|$.  It follows from (4.16) of \cite{Tadmor} that
    $$\|\bx-\bp\|_{\bH^1}^2 \le \sum_{j=1}^d \left(8 C_\lam |\bx(-1)| (\sinh\eta)^{-1} k e^{-k\eta}\right)^2 = d \left(8 C_\lam (\sinh\eta)^{-1} k e^{-k\eta}\right)^2 |\bx(-1)|^2.$$
Because $|\bx(-1)|\le 0.9062 \|\bx\|_{\bH^1} = 0.9062$, we are done.\end{proof}

\subsection{Collocation approximation of the monodromy operator and its eigenvalues}  Define
    $$\hM_B = \begin{pmatrix}B(t_0) & & & \\ & \ddots & & \\ & & B(t_{N-1}) & \\ I_d & & & 0_d \end{pmatrix}$$
where $I_d,0_d$ are the $d\times d$ identity and zero matrices, respectively.  The approximate monodromy matrix for DDE \eqref{DDEsys} is
    $$\bU_N = \left(\hbD_N - \hM_A\right)^{-1} \hM_B \in \CC^{l\times l},$$
where $l=d(N+1)$, assuming the inverse exists (which we determine by numerical means).

Let $\bPi_N:\bH^1\to\Po_N$ be the orthogonal projection onto $\Po_N$, so that if $\bbf(t)=\sum_{j=0}^\infty \balpha_j \tilde T_j(t)$ then $\bPi_N \bbf = \sum_{j=0}^N \balpha_j\tilde T_j(t)$.  We define the approximate monodromy operator
    $$\hbU_N = \bP_N \bU_N \bE_N \bPi_N \in \Lop(\bH^1).$$
The following lemma is the systems generalization of lemma \ref{lem:diagUm} in subsection \ref{subsect:aposteig}.

\begin{lem}\label{lem:diagUNsys}  Let $l=d(N+1)$.  Suppose $\bU_N = \bV \bLambda \bV^{-1}$ with $\bLambda=(\lam_j)_{j=1}^l$ diagonal.  Let $\bv_k$, $k=1,\dots,l$, be the $k$th column of $\bV$ and define $\bp_k=\bP_N \bv_k$, a $\CC^d$-valued polynomial of degree $N$.

Expand $\bp_k$ in Chebyshev series $\bp_k(t) = \sum_{j=1}^{N+1} \sum_{s=1}^d \Gamma_{(j-1)d+s,k} \tT_{j-1}(t)\otimes \be_s$ where $\Gamma_{i,k} \in \CC$ for $i,k=1,\dots,l$.  The matrix $\bGamma=\left(\Gamma_{i,k}\right)_{i,k=1}^l$ is invertible and, in fact, if
    $$\bC =\sqrt{\pi/2} \left(\operatorname{diag}(2^{-1/2},2,3,\dots,N+1) \,C\right)\otimes I_d,$$
where $C$ is the matrix described by equation \eqref{discreteCheb}, then $\bGamma=\bC \bV$.

Let $\bl^2 = l^2\otimes \CC^d$ and denote a typical element by $\ba=(a_1^1,\dots,a_1^d,a_2^1,\dots,a_2^d,\dots)$, $a_j^s \in \CC$.  Let $\hbV:\bl^2\to\bH^1$ be defined by
\begin{align*}
\hbV\ba &= \sum_{j=1}^{N+1} \sum_{s=1}^d a_j^s \bp_{(j-1)d+s} + \sum_{j>N+1} \sum_{s=1}^d a_j^s \tT_{j-1}(t) \be_s \\ 
    &= \sum_{j=1}^{N+1} \sum_{s=1}^d \sum_{j'=1}^{N+1} \sum_{s'=1}^d a_j^s \Gamma_{(j'-1)d+s',(j-1)d+s} \tT_{j'-1}(t) \be_{s'} + \sum_{j>N+1} \sum_{s=1}^d a_j^s \tT_{j-1}(t) \be_s.
\end{align*}
Then $\hbV$ is bounded and invertible and furthermore
    $$\|\hbV\| \le \left(\|\bGamma\|_2^2 + 1\right)^{1/2}, \qquad \|\hbV^{-1}\| \le \left(\|(\bGamma)^{-1}\|_2^2 + 1\right)^{1/2}.$$

Finally, define $\hbLambda\in\Lop(\bl^2)$ by $(\hbLambda \ba)_{(j-1)d+s} = \lam_{(j-1)d+s} a_j^s$ if $j<N+1$ and $s=1,\dots,d$.  If $k>l=d(N+1)$ then let $(\hbLambda \ba)_{k} = 0$.  (Thus $\hbLambda$ is a diagonal operator on $\bl^2$ of rank at most $l$.)  Then $\hbU_N = \hbV \hbLambda \hbV^{-1}$ and $\|\hbU_N\| \le \left(\max_{1\le j \le l} |\lam_j|\right)\,\cond(\hbV)$.\end{lem}

\begin{proof} Follow the proof of lemma \ref{lem:diagUm} with appropriate modifications. \end{proof} 

We now give our main theorem on the monodromy operator $\bU$ for DDE \eqref{DDEsys}.

\begin{thm}\label{thm:mainsys}  Let $d\ge 1$ and $N\ge 1$.  Suppose $A,B$ in DDE \eqref{DDEsys} are analytic $d\times d$ matrix-valued functions with common regularity ellipse $E\supset I$.  Let $e^\eta = S+s >1$ be the sum of semiaxes of $E$.  Suppose $C_\lam$ is an upper bound for $|\Phi_\lam(z)|$ for $z\in E$ and $\Phi_\lam(z)$ the analytic continuation of $\Phi_\lam(t)$, where $\Phi_\lam(t)$ is the fundamental solution of $\dot\bx = (A(t)+B(t)/\lam) \bx(t)$.

Let $\delta>0$.  Suppose $\hbU \bx = \mu \bx$ for $\bx\in \bH^1$, $\|\bx\|=1$, and $|\mu|\ge \delta$.  Let
\begin{equation}\label{epsksys}
\eps_k = 8 \sqrt{d} C_\lam (\sinh\lam)^{-1} k\,e^{-k\eta}
\end{equation}
for $k=1,\dots,N$.  Assume $\bU_N = \bV\bLambda \bV^{-1}$ with $\bLambda = \left(\lam_i\right)_{i=1}^l \in \Lop(\bl^2)$ diagonal and eigenvalues $\lam_i$ ordered $|\lam_1|\ge |\lam_2|\ge \dots$.  As in lemma \ref{lem:diagUNsys} it follows that $\hbU_N = \hbV \hbLambda \hbV^{-1}$.

Suppose $\|\hbU(\tT_{j-1}\otimes\be_s) - \hbU_N(\tT_{j-1}\otimes\be_s)\|_{\bH^1} \le \nu_j^s$ where $j=0,\dots,N$ and $s=1,\dots,d$.  Let $\xi_k^2 = \sum_{j=0}^k \sum_{s=1}^d (\nu_j^s)^2$ and 
    $$\omega_k = \eps_k \left(\|\hbU\|+|\lam_1|\cond(\hbV)\right) + (1+\eps_k)\xi_k$$
for $k=1,\dots,N$.  Then
\begin{equation}\label{thmresultsys}
\min_{i=1,\dots,d(N+1)} |\mu-\lam_i| \le \min\left\{\omega_1,\dots,\omega_N\right\} \cond(\hbV).\end{equation}
\end{thm}

\begin{proof}  By lemma \ref{lem:eigsys}, for each $k=1,\dots,N$ we have $\|\bx-\bq_k\|_{\bH^1} < \eps_k$ where $\bq_k=\dot \bI_k \bx_k$.  Apply corollary \ref{cor:BFHwp} with $\Xcal=\bH^1$, $A=\hbU_N$, $B=\hbU$, $x=\bx$, $p=\bq_k$, and $\eps=\eps_k$.  Conclude
\begin{equation}\label{firstconcludesys}
\min_{i=1,\dots,d(N+1)} |\mu-\lam_i| \le \cond(\hbV) \left(\eps_k(\|\hbU\|+\|\hbU_N\|) + \|(\hbU-\hbU_N)\bq_k\|\right).\end{equation}

On the other hand, $\bq_k=\sum_{j=0}^k \sum_{s=1}^d \alpha_{kj}^s \tT_j(t)\otimes\be_s$ so
    $$\|(\hbU-\hbU_N)\bq_k\|_{\bH^1} \le \sum_{j=0}^k \sum_{s=1}^d |\alpha_{kj}^s| \nu_j^s \le \|\bq_k\|_{\bH^1} \xi_k \le \left(\|\bx\|_{\bH^1} + \|\bq_k-\bx_k\|_{\bH^1}\right) \xi_k \le (1+\eps_k)\xi_k$$
by Cauchy-Schwarz.  From \eqref{firstconcludesys} and lemma \ref{lem:diagUNsys} we conclude with estimate \eqref{thmresultsys}.\end{proof}

\begin{rem}  As noted is section \ref{subsect:collocalgsys} we can produce good \emph{a posteriori} (``bootstrapping'') bounds $C_A$ on $|\Phi_A(t)|$ and $|\Phi_A(t)\Phi_A(s)^{-1}|$ where $\Phi_A(t)$ is the fundamental solution to $\dot\bx=A(t)\bx$.  This has two consequences for applications of theorem \ref{thm:mainsys}.  First, the \emph{a posteriori} quantities $\nu_j^s$ will be small (comparable to machine precision in practice; see theorem \ref{thm:apostsys}).  Second, the estimate for $\|\hbU\|$ will be reasonable (see lemma \ref{lem:Uhatnormsys}).

On the other hand, the bound $C_\lam$ on the analytic continuation $\Phi_\lam(t)$ of the fundamental solution $\Phi_\lam(t)$ to $\dot\bx = \left(A(t)+B(t)/\lam\right) \bx$ will still be \emph{a priori}.\footnote{One might imagine seeking to find $\Phi_\lam(z)$ for $z\in E$ by a spectral method.  Unfortunately, this problem is a system of first order \emph{partial} differential equations.  We don't need an accurate solution to it, in any case.  Perhaps an \emph{a posteriori} bound can be sought by some practical method; we don't know of one yet.}  Such a bound will inevitably be big (see the example which follows).  Thus $\eps_k$ will be big for small $k$.  Nonetheless the essential point of figure \ref{fig:constantsdecay} remains: the consistently small values for $\nu_j^s$ will allow us to choose $k\le N$ large enough so that $\eps_k$ is small.  (One must choose $N$ large enough, of course.  This can be done using equation \eqref{epsksys} for $\eps_k$.)

Said another way, good \emph{a posteriori} estimates on initial value problems eventually wins when competing with large ignorance of the regularity of eigenfunctions of $\hbU$.\end{rem}

\section{Example: a delayed, damped Mathieu equation}  \label{sect:Mathieu}

\begin{figure}[ht]
\regfigure{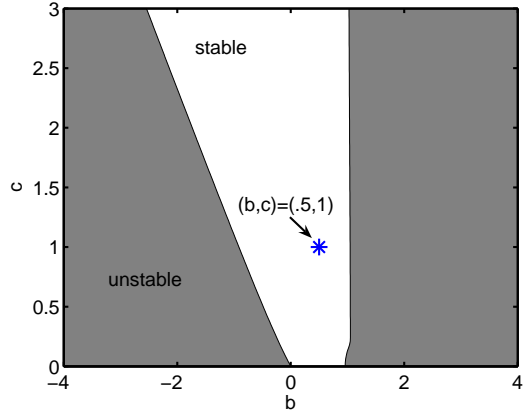}{2.3}
\caption{Stability chart for the delayed, damped Mathieu equation $\ddot x + c \dot x + (1 + \cos(\pi t))x=b x(t-2)$.} \label{fig:Mathieuchart}
\end{figure}

The DDE $\ddot x + c \dot x + (1 + \cos(\pi t))x=b x(t-2)$, $b\in\RR$, $c\ge 0$ is a delayed, damped Mathieu equation (compare \cite{InspergerStepan02B}).  It corresponds to $\dot\by = A(t) \by + B(t) \by(t-2)$ where 
    $$A(t)=\begin{pmatrix} 0 & 1 \\ -1-\cos(\pi t) & -c\end{pmatrix} \qquad \text{and} \qquad B(t)=\begin{pmatrix} 0 & 0 \\ b & 0\end{pmatrix}.$$
For $(b,c) \in [-4,4]\times [0,3]$ we have the numerically-produced stability chart of figure \ref{fig:Mathieuchart} and we expect that the particular parameter pair $(b,c)=(1/2,1)$ will be stable.  In fact, for this parameter pair we get the numerical value $|\mu_1|=0.5858$ for the largest eigenvalue of $\bU$.  The machinery of the previous section allows us to prove the correctness of this eigenvalue to two decimal places, in an \emph{a posteriori} manner, as follows.

First, using the same method as example \ref{example:getboundMathieu}, we find an \emph{a priori} bound $|\Phi_\lam(z)| \le C_A = e^{2\sqrt{6}}\approx 134.2$ on the fundamental solution $\Phi_A(t)$ solving $\dot \Phi_A = A(t) \Phi_A$, $\Phi_A(-1)=I$ on $t\in[-1,1]$.  We improve this bound by \emph{a posteriori} iterations as in subsection \ref{subsect:collocalgsys} to a new bound $C_A' = 5.12$.

If $\bx$ is a ($\CC^2$-valued) eigenfunction of $\bU$ corresponding to eigenvalue $\lam$ for $|\lam|\ge \delta >0$ then $\dot \bx = (A+B/\lam)\bx$ and $\|A+B/\lam\|_\infty \le \sqrt{1+c^2+(2+|b|/\delta)^2}$ using the Frobenius norm.  We need to bound the analytic continuation $\Phi_\lam(z)$, for $z$ in a regularity ellipse $E$, of the fundamental solution $\Phi_\lam(t)$ of $\dot \bx = (A+B/\lam)\bx$.  Let $s=0.5$ and $S=\sqrt{1+s^2}=\sqrt{1.25}$ be semiaxes of $E$ where $E$ is an ellipse with foci $\pm 1$.  Then $|\Phi_\lam(z)|$ is bounded by an \emph{a priori} estimate $|\Phi_\lam(z)| \le \exp((1+S) \|A+B/\lam\|_\infty) \le 4121 = C_\lam$ for $z\in E$ if $\delta=0.3$, $b=0.5$, and $c=1$.

\begin{figure}[ht]
\regfigure{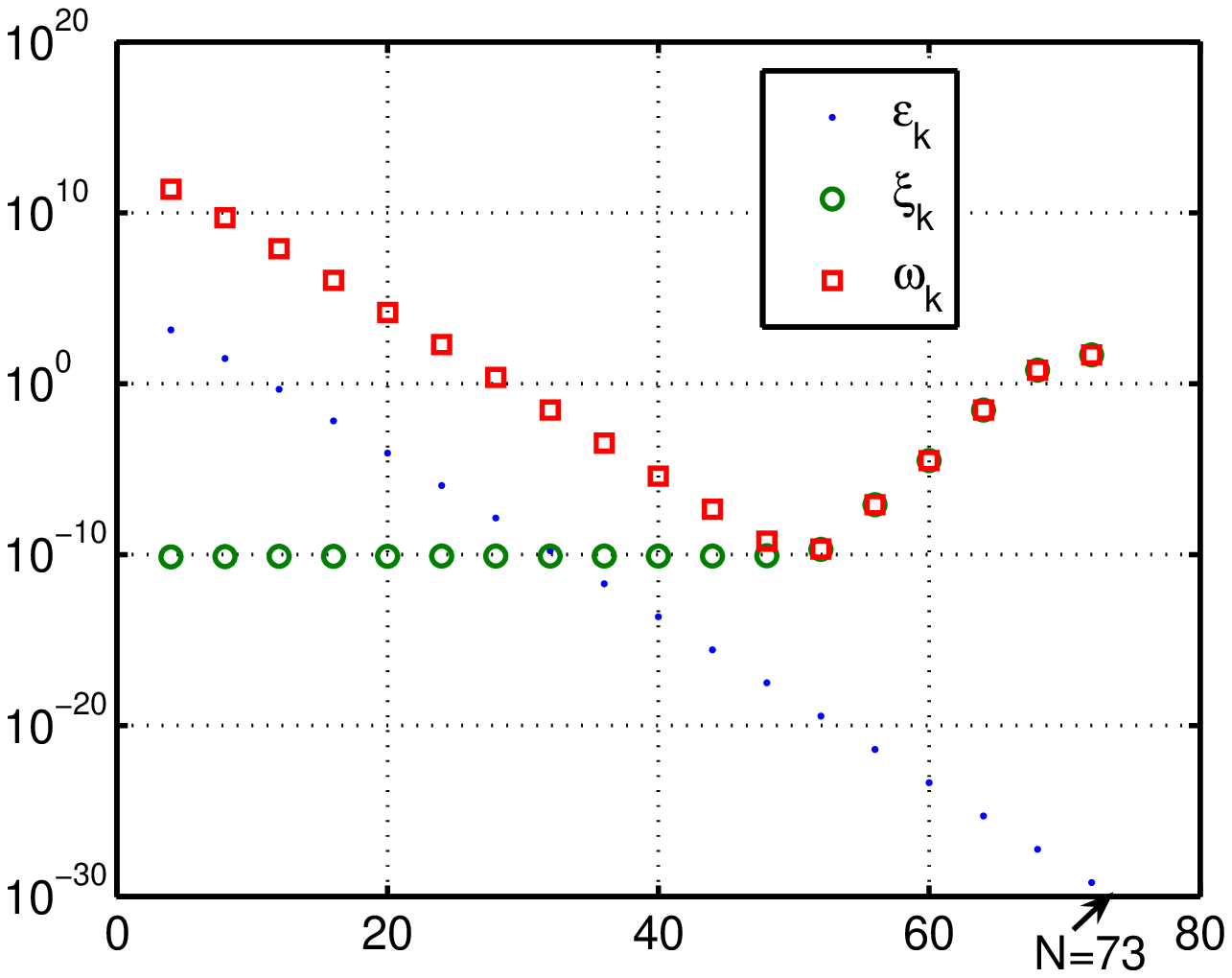}{2.35}\regfigure{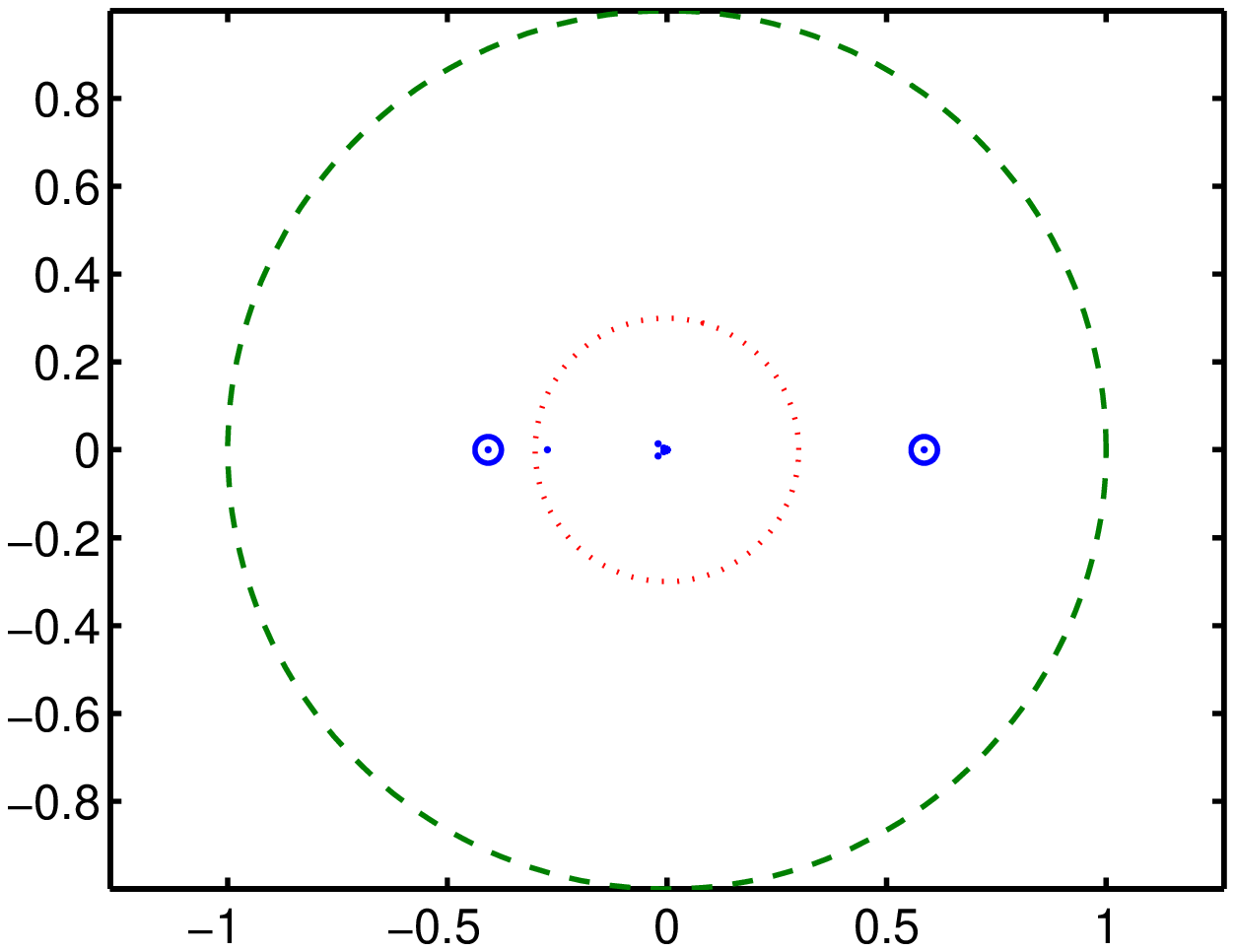}{2.35}
\caption{\textbf{(a)}  Behavior of constants in applying theorem \ref{thm:mainsys} to our delayed, damped Mathieu equation.  \textbf{(b)}  Two (numerical) eigenvalues of $\bU_N$, $N=73$, of our delayed, damped Mathieu equation exceed $\delta=0.3$, as shown.  Theorem \ref{thm:mainsys} proves that any eigenvalues of the monodromy operator $\bU$ will lie in the discs of radius $0.03019$ around the numerical eigenvalues, as shown.} \label{fig:mathieueigs}
\end{figure}

We now apply theorem \ref{thm:mainsys} with (after trying some smaller values) $N=73$.  The constants $\eps_k$, $\xi_k$, and $\omega_k$ which result are shown in figure \ref{fig:mathieueigs}\textbf{a}.  The condition number estimate for $\hbV$ given in lemma \ref{lem:diagUNsys} is computed to be approximately $2.898 \times 10^8$, and this limits our expectations for accuracy of the numerical eigenvalues to perhaps 6 or 7 digits (given that double precision is about 15 digits).  The error radius (the right side of equation \eqref{thmresultsys} and the main result of theorem \ref{thm:mainsys}) is $0.03019$.  We have the picture of eigenvalues in figure \ref{fig:mathieueigs}\textbf{b}, complete with error bounds.  The eigenvalues of the exact monodromy operator $\bU$ which exceed $\delta=0.3$ in magnitude (shown) are proven to lie within the discs of radius $0.0319$ (shown) around the numerically computed eigenvalues (shown).

We have proven the stability of this delayed, damped Mathieu equation for the parameter pair $(b,c)=(0.5,1)$.

\section{Techniques for evaluation of polynomials and norms}\label{sect:polytechs}

Consider evaluating a polynomial $p\in\Po_N$ described by its values $v\in \CC^{N+1}$ at $N+1$ Chebyshev collocation points.  Concretely, suppose we want to find $z=p(t)$ for $t\in [-1,1]$ given $v_j=p(t_j)$ at the points $t_j=\cos(j\pi/N)$.

An easy method is described by the single line of \Matlab:\begin{quote}\texttt{
z=polyval(polyfit(tj,v,N),t)
}\end{quote}
This method produces considerable error for large $N$ because \texttt{polyfit} finds the coefficients of $p$ in the monomial basis, an ill-conditioned operation on the values $v$ \cite{Trefethen}.   

A good alternative is the ``barycentric'' interpolation of \cite{Salzer} (see also \cite{BattlesTrefethen}).  The following simple formula is special to Chebyshev collocation (extreme) points:
    $$p(t)=\sum_{j=0}^N \frac{w_j}{t-t_j} v_j \bigg/ \sum_{j=0}^N \frac{w_j}{t-t_j}, \quad \text{where} \quad w_j=\begin{cases} (-1)^j/2, & j=0 \text{ or } j=N,\\ (-1)^j, &\text{otherwise.}\end{cases}$$
A proof of the numerical stability of this method appears in \cite{Higham}; see also \cite{BerrutTrefethen}.  The following is a brief \Matlab~implementation, though not necessarily an optimal one.
\small\begin{quote} \begin{alltt} 
function p=bary(v,t,N)
%BARY Barycentric interpolation for Chebyshev points.  
tj=cos(pi*(0:N)/N)';                                     % Cheb points
wj = [1/2; ones(N-1,1); 1/2].*(-1).^((0:N)');
for k=1:length(t)
   if min(abs(t(k)-tj)) > 2*eps                          % at generic pt
      zz=wj./(t(k)-tj); p(k)=sum(zz.*v)/sum(zz);
   else, j=find(abs(t(k)-tj) <= 2*eps); p(k)=v(j); end   % very near Cheb pt
end
\end{alltt} \end{quote}\normalsize\smallskip

\begin{figure}[ht]
\regfigure{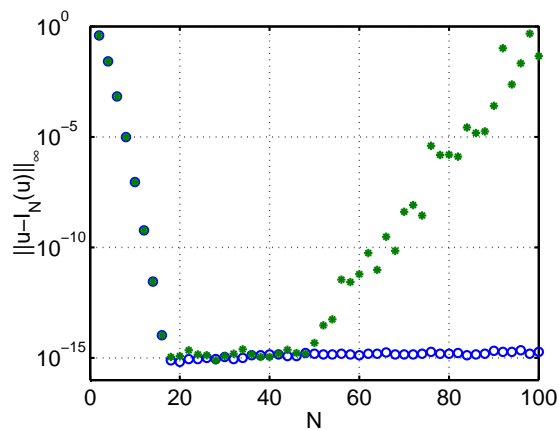}{2.3}
\caption{Importance of the right polynomial evaluation method: Approximations of $\|u-I_N u\|_\infty$ where $u(t)=\sin(2t)$: dots are by ``\texttt{polyval(polyfit(\dots))}'' and circles are by \texttt{bary.m}, above.} \label{fig:polybary}
\end{figure}

\begin{example}  To estimate the norm $\|u-I_N u\|_\infty$, one may evaluate $|u(t)-p(t)|$ at many equally-spaced points $t$ in $[-1,1]$ where $p=I_N u$.  Specifically, let  $u(t)=\sin(2t)$.  Figure \ref{fig:polybary} shows that the evalution of $p(t)$ is done stably by barycentric interpolation up to degree $N=100$, but that the ``\texttt{polyval(polyfit(\dots))}'' method behaves catastrophically for large $N$.\end{example}

The above is really just a warm-up.  The polynomial calculation most used in this paper is the estimation of $\|p\|_\infty$ when $p(t)$ is a polynomial represented by its collocation values, or of $\|\alpha\|_\infty$ when $\alpha(t)$ is an analytic function.  For estimating such norms there is a much better method than evaluating $p(t)$, respectively $\alpha(t)$, at predetermined points.

Consider the polynomial case.  If $p(t)=\sum_{k=0}^N a_k T_k(t)$ where $T_k(t)=\cos(k\arccos t)$ is the $k$th Chebyshev polynomial then
\begin{equation}\label{coeffest}
\|p\|_\infty \le \sum_{k=0}^N |a_k|
\end{equation} 
since $|T_k(t)|\le 1$.  Estimate \eqref{coeffest} requires the coefficients $a_k$ in the Chebyshev expansion of $p$ but these are easily calculated by the FFT, as for instance by the \Matlab~function 
\small \begin{quote} \begin{alltt} 
function a = coefft(v)
%COEFFT  Compute Cheb coefficients of an N degree polynomial p(t) 
%represented by its colloc values v_j = p(t_j) where t_j = cos(pi j/N)
N = length(v)-1; if N==0, a=0; return, end
v = v(:);                                    % force into column
U = fft([v; flipud(v(2:N))])/N;              % do t -> theta then FFT
a = ([.5 ones(1,N-1) .5])'.*U(1:N+1);
\end{alltt} \end{quote}\smallskip\normalsize
Thus $\|p\|_\infty \le$ \texttt{sum(abs(coefft(v)))} if \texttt{v} is a vector of the collocation values.

If $\alpha(t)$ is analytic on $I$ then we estimate $\|\alpha\|_\infty$ by approximating $\alpha(t)$ by high degree Chebyshev interpolation and then we use \eqref{coeffest}.  A very similar issue is addressed in \cite{BattlesTrefethen}; we resolve ours in a similar way, as follows.

Given $\alpha(t)$, we start with $M=15$ (for concreteness), evaluate at collocation points $t_j=\cos(\pi j/M)$, $j=0,\dots,M$, and use the FFT to calculate the coefficients $a_k$ of the corresponding polynomial.  We then inquire if the last few coefficients (say, the last four, for concreteness) are small.  In particular, if $\max\{|a_{M-3}|,\dots,|a_M|\} < 10 \eps_m$, where $\eps_m$ is machine precision, then we accept the estimate $\|\alpha\|_\infty \approx \|p\|_\infty$ and use estimate \eqref{coeffest}.  If not, $M$ is doubled (actually, to one less than the next power of $2$; this gives efficiency in the FFT) and we try again.  We stop with $M=2^{12}-1$.

\begin{example}  Again consider estimating $\|u-I_N u\|_\infty$, this time with $N=5$ fixed.  Thus $\alpha=u - I_5 u$ is the analytic function in question.  The following \Matlab~implements the above algorithm:
\small \begin{quote} \begin{alltt} 
b = coefft(u(cos(pi*(0:5)/5)));                % coeffs of I_N u
M = 15;  last = 8;  tol = 10*eps;
for s = 1:last
    tM = cos(pi*(0:M)/M);
    a = coefft(u(tM)) - [b; zeros(M-5,1)];     % coeffs of u - I_N u
    z = sum(abs(a));
    if max(abs(a(M-2:M+1))) < tol, break, end
    M = 2*(M+1)-1
end
\end{alltt}  \end{quote}\smallskip\normalsize
If $u(t)=\sin(2t)$, this procedure stops at $N=31$ with an estimate $\|u-I_5 u\|_\infty \le 0.00070975$.  By contrast, the use of 1000 equally-spaced points in $[-1,1]$ and \texttt{bary.m}, as in
\small\begin{quote} \begin{alltt} 
tp=linspace(-1,1,1000); up=sin(2*tp);
tj=cos(pi*(0:5)/5)'; v=sin(2*tj);
for k=1:1000, INup(k)=bary(v,tp(k),5); end
\end{alltt}  \end{quote}\smallskip\normalsize
yields $\|u-I_N u\|_\infty \approx$ \texttt{max(abs(up-INup))} $= 0.00067538$ at substantially greater cost (roughly 100 times the execution time). \end{example}

\section{The monodromy operator for the $(\text{delay}) \ne (\text{period})$ case}\label{sect:taunoteqperiod}

As noted in the introduction and on page \pageref{easycasepage}, if the period is equal to (or an integer multiple of) the delay, \emph{and} if the Floquet matrix of a corresponding homogeneous ODE is exactly known,\footnote{The Floquet matrix for a linear system of non-constant-coefficient ODEs is usually not exactly known.  Thus the bulk of the paper, and section \ref{sect:systems} in particular, addresses a nontrivial case.} then complex variable methods can be applied to the stability problem for linear, periodic, fixed-delay DDEs; see section 8.3 of \cite{HaleLunel}.  In this section we write down a variation-of-parameters formula for the monodromy operator in (essentially) the general case of one fixed delay and periodic coefficients.  We assume no rational relation between period and delay.  We see that our methods generalize.  The same \emph{a posteriori} results for numerical initial value problems and for numerical approximations of the eigenvalues of the monodromy operators should be available with sufficient elaboration.

Consider the system of DDEs (compare subsection \ref{subsect:DDEsystemmono}):
\begin{equation}\label{DDEsysdelaynotperiod} 
\dot\by(t) = A(t)\by(t) + B(t) \by(t-\tau)
\end{equation}
for $\by(t)\in\CC^d$ and $A,B$ continuous matrix-valued (i.e.~in $\CC^{d\times d}$) functions with period $T$ \emph{exceeding} $\tau$ \emph{but less than or equal to} $2\tau$:
    $$0<\tau < T\le 2\tau, \quad A(t+T)=A(t), \quad B(t+T)=B(t).$$
Other cases than $\tau < T \le 2\tau$ will be handled by appropriate modification of the argument below.  Note $T$ can always be forced to exceed $\tau$ by choosing the nominal coefficient period to be a sufficiently large integer multiple of the minimal period.

Now suppose $\bbf(t)$ is a history function defined on the interval $-\tau \le t \le 0$.  The solution of \eqref{DDEsysdelaynotperiod} for one full period $T$, with initial value $\by(t)=\bbf(t)$, $t\in [-\tau,0]$, is found by solving two consecutive ODE problems
\begin{align}
    &\dot \by_1 = A(t) \,\by_1 + B(t) \,\bbf(t-\tau), &&\by_1(0)=\bbf(-\tau), &&t\in [0,\tau], \label{twoivps} \\
    &\dot \by_2 = A(t) \,\by_2 + B(t) \,\by_1(t-\tau), &&\by_2(\tau)=\by_1(\tau), &&t\in [\tau,T]. \notag
\end{align}
because $2\tau$ exceeds $T$.

\begin{figure}[ht]
\widefigure{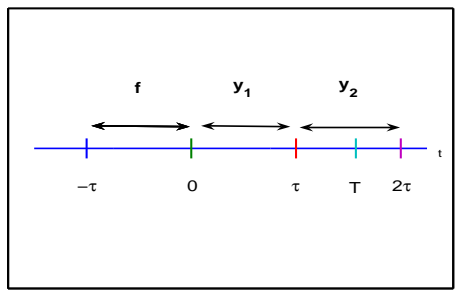}{1.8}{3.0}
\caption{The initial value problem for \eqref{DDEsysdelaynotperiod} in the case $0<\tau<T\le 2\tau$.} \label{fig:twoivps}
\end{figure}

In figure \ref{fig:twoivps} we show the schematic map from $\bbf$ to the solution $\by_1,\by_2$ on $[0,T]$.  However, to be considered as the action of an operator possessing eigenvalues, the solution of the initial value problem for system \eqref{DDEsysdelaynotperiod} must be redefined to become an operator acting from a space to the same space.

Let $\Ycal = \bC[0,T-\tau]\oplus \bC[T-\tau,T]$, that is, let $\Ycal$ be $\bC[0,T]$ broken at $t=T-\tau$.  We denote a typical element of $\Ycal$ by $(\bg,\bbf)$ with $\bg\in \bC[0,T-\tau]$ and $\bbf\in \bC[T-\tau,T]$.

\begin{figure}[ht]
\widefigure{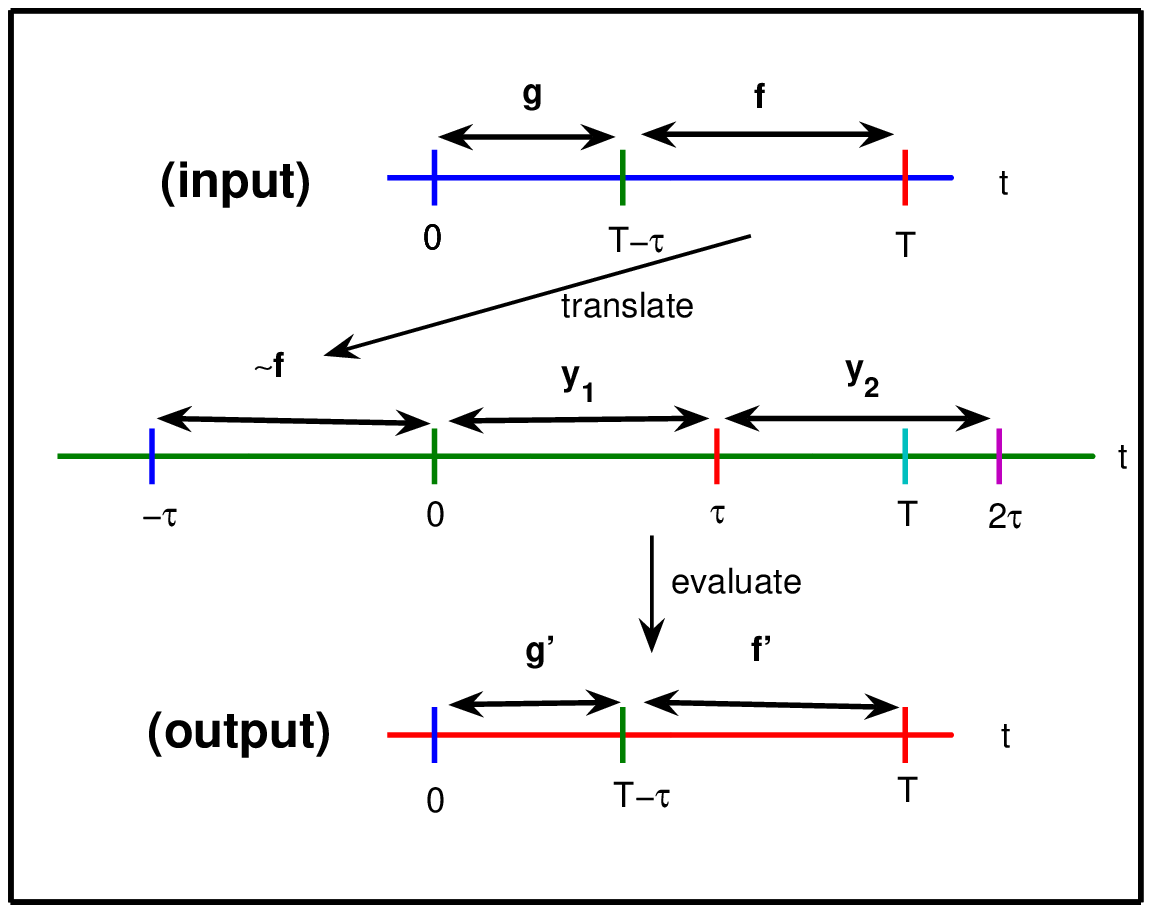}{2.5}{3.5}
\caption{Schematic of action of $\bU$ on $\Ycal$.} \label{fig:schemeforU}
\end{figure}

The monodromy operator can now be defined as the following composition
\begin{align*}
\bU &: (\bg,\bbf)\in\Ycal \overset{\text{ignor } \bg}{\mapsto}  \bbf\in \bC[0,T-\tau]  \overset{\text{translate}}{\mapsto}  \tilde\bbf\in \bC[-\tau,0]  \\
    &\qquad \overset{\text{solve IVP by \eqref{twoivps}}}{\mapsto}  \left(\by_1\in \bC[0,\tau],\by_2\in \bC[\tau,2\tau]\right)  \overset{\text{evaluate onto } [0,T]}{\mapsto}  (\bg',\bbf')\in\Ycal.
\end{align*}
This composition is pictured in figure \ref{fig:schemeforU}.

The ``evaluate onto $[0,T]$'' stage perhaps requires explanation.  We have
    $$\bg'(t)=\by_1(t), \quad t\in[0,T-\tau], \qquad \text{and} \qquad \bbf'(t) = \begin{cases} \by_1(t), & t\in[T-\tau,\tau], \\ \by_2(t), & t\in[\tau,T].\end{cases}$$

Note that though $(\bg,\bbf)$ need not be in $\bC[0,T]$ (i.e.~$g(T-\tau)$ need not equal $\bbf(T-\tau)$), it is still true that $\by_2(\tau)=\by_1(\tau)$ and thus we may say $(\bg',\bbf') \in \bC[0,T]$.  In particular, the eigenfunctions of $\bU$ will be continuous functions.

Furthermore, because the input $\bg$ is ignored, $\bU$ has an infinite-dimensional null space and thus $0\in\sigma(\bU)$ has infinite multiplicity.  Nonetheless, $\bU$ is evidently compact.  Powers of $\bU$ solve initial value problems for equation \eqref{DDEsysdelaynotperiod} out to arbitrary $t>0$.  The nonzero eigenvalues of $\bU$ determine the stability of \eqref{DDEsysdelaynotperiod}.

We propose that a numerical approximation of $\bU$ can be constructed by using separate polynomial approximations of the inputs $\bg$ and $\bbf$ by interpolation at the (scaled and shifted) Chebyshev collocation points in intervals $[0,T-\tau]$ and $[T-\tau,\tau]$, respectively.  (In fact, we need do no actual interpolation of $\bg$.)  The ODE initial value problems in \eqref{twoivps} are then solved by the method of section \ref{sect:apost}, so that \emph{a posteriori} estimates are available for the approximations $\bp_1,\bp_2$ to $\by_1,\by_2$.  Note $\bp_1,\bp_2$ are functions on $[0,\tau]$ and $[\tau,2\tau]$, respectively.  We then interpolate to find $\bq_1\approx \bg'$ on $[0,T-\tau]$ and $\bq_2\approx\bbf'$ on $[T-\tau,T]$.  The information in $\bp_2$ on $[T,2\tau]$ is lost by this procedure, a small loss in efficiency if the goal is a long-time solution of the initial value problem.  (There is no loss of accuracy, however, \emph{if} $A,B,\bbf$ are smooth on $(-\eps,T+\eps)$, $(-\eps,T+\eps)$, and $(-\tau-\eps,+\eps)$, respectively, because $\by_1,\by_2$ connect smoothly at $t=\tau$ in this case.)  The result of these operations is a numerical approximation to $\bU$ acting from the space $\Pcal_{N_1}\oplus \Pcal_{N_2}$ of polynomials of (possibly) different degrees $N_1,N_2$ back to itself.

Though we admit that ``bookkeeping'' complications have kept us from filling out the details in these limited pages, the \emph{a posteriori} and eigenvalue perturbation analysis of earlier sections will first give estimates of the accuracy of initial value problems for DDE \eqref{DDEsysdelaynotperiod} and will furthermore give eigenvalue estimates as well.

\section{Conclusion}\label{sect:conclude}

We conclude with a sketch of the ``roads not taken'' in the current paper and also with certain genuinely open problems.

First, this paper does not address the \emph{multiple} fixed delays case for linear, periodic DDEs.  This omission is merely a technical one; all techniques here generalize without difficulty to that case.

Second, we have not yet considered the case of linear, periodic DDEs with nonconstant (but state-independent) delays (i.e.~$\dot \bx(t)=A(t) \bx(t) + B(t) \bx(t-\tau(t))$), or the case of neutral linear, periodic DDEs ($M_0(t) \dot \bx(t) + M_1(t) \dot \bx(t-\tau)=A(t) \bx(t) + B(t) \bx(t-\tau)$), or the case of linear, periodic functional equations ($\dot \bx(t) = \int_0^{\tau(t)} k(t,s) \bx(t-s) \,d\nu(s)$; \cite{HaleLunel}).  It is clear that many techniques in the current paper will generalize if we assume analyticity for the relevant functions ($\tau(\cdot)$, $M_i(\cdot)$, and $k(\cdot,s)$, respectively).  We make this statement confidently but without precision.

Now, specific open problems:\begin{itemize}
\item Under what conditions does $\bU$ or $\hbU$ (the exact monodromy operator, in any either case) actually diagonalize?  
\item In fact, do the eigenfunctions of $\hbU$ generically form a  Riesz basis \cite{GohbergKrein} or some other (substantially) linearly-independent basis for $\bH^1$?  (Better yet, find an \emph{a priori} estimate of the condition number of the finite rank operator with columns the normalized eigenfunctions of $\hbU$ corresponding to large eigenvalues $\mu$, e.g.~ such that $|\mu|\ge \delta$ for fixed $\delta\ge 0$.)
\item Given $N$, estimate $\|\hbU-\hbU_N\|_{\bH^1}$.
\end{itemize}

For none of these questions would an answer actually effect the \emph{a posteriori} methods of the current paper.  Indeed these questions relate to building an \emph{a priori} understanding of a spectral approximation of monodromy operators for linear DDEs.


\bibliography{../ddec}
\bibliographystyle{siam}

\small
\medskip

\phantom{bob}

\addcontentsline{toc}{section}{Appendices}

\centerline{\textsc{Appendix A: Constants for theorem \ref{thm:apostivp} when $a(t)=a_0$}}
\medskip

To prove the special-case estimates \eqref{specialone} and \eqref{specialtwo}, we exploit cancellation in the second integral appearing in \eqref{inttoestimate}.  We need a lemma which is well-known in the theory of Fourier series, namely, that integration-by-parts allows us to estimate the cancellation which occurs in an integral of a regular function $f(x)$ against a rapid sinusoid.  Actually, we will need the given corollary.

\begin{lem*}\label{lem:cancel}  If $f\in C^1[a,b]$ then
    $$\left|\int_a^b f(x)\begin{Bmatrix} \sin(nx) \\ \cos(nx) \end{Bmatrix}\,dx \right| \le \frac{1}{n} \left(2 \|f\|_\infty + (b-a) \|f'\|_\infty\right).$$
\end{lem*}

\begin{proof}  Integrate-by-parts.  \end{proof}

\begin{cor*}\label{cor:lNcancel}  If $t\in I$ then
    $$\left|\int_{-1}^t e^{a_0(t-s)} l_N(s)\,ds\right| \le \frac{(\pi(|a_0|+1)+4) e^{2|\Re a_0|}}{N 2^{N-1}}.$$
\end{cor*}

\begin{proof}  Substitute $s=\cos\theta$ and use lemma \ref{lem:lNtotrig} to find
\begin{align}\label{inttocancel}
\int_{-1}^t e^{a_0(t-s)} l_N(s)\,ds &= \frac{-e^{a_0 t}}{2^{N-1}} \int_{\arccos t}^\pi e^{-a_0 \cos\theta} \left(1-\cos(\theta)\right) \sin(N\theta)\,d\theta.
\end{align}
Let $f(\theta)=e^{-a_0 \cos\theta} \left(1-\cos(\theta)\right)$ so $f\ge 0$,  $\|f\|_\infty\le 2 e^{|\Re a_0|}$, and $\|f'\|_\infty \le (|a_0|+1) e^{|\Re a_0|}$.  By lemma \ref{lem:cancel},
    $$\left|\int_{\arccos t}^\pi f(\theta) \sin(N\theta)\,d\theta\right| \le  \frac{1}{N} (\pi(|a_0|+1)+4)e^{|\Re a_0|}$$
and the result follows.\end{proof}

If $\Re a_0\le 0$ and $N$ is small then we can sometimes get a slightly better estimate by not pursuing cancellation.  In fact, if $\Re a_0\le 0$ then
\begin{align*}
\left|\int_{-1}^t e^{a_0(t-s)} l_N(s)\,ds\right| &\le \int_{-1}^t |l_N(s)|\,ds = \frac{1}{2^{N-1}} \int_{\arccos t}^\pi \left|\frac{(\cos(\theta)-1) \sin(N\theta)}{\sin(\theta)}\right| \sin\theta\, d\theta \\
&\le \frac{1}{2^{N-1}} \int_0^\pi 1-\cos\theta\,d\theta = \frac{\pi}{2^{N-1}}
\end{align*}
by lemma \ref{lem:lNtotrig}.

We can now start with equation \eqref{inttoestimate} and complete the proof of theorem \ref{thm:apostivp} in this special case.  Since $I_N(ap)=a_0 p$, and because $|\Phi_s^t|=|e^{a_0(t-s)}|\le \max\{e^{2 \Re a_0},1\}$ for $-1\le s\le t\le 1$, it follows that
    $$\|y-p\|\le c_1 \|u-w\|_\infty + |\beta| \max_{t\in I} \left|\int_{-1}^t e^{-a_0(t-s)} l_N(s)\,ds\right|$$
where $c_1 = \max\{2 e^{\Re a_0},2\}$.  By corollary \ref{cor:lNcancel} and equation \eqref{initialfact}, $\|y-p\|_\infty\le c_1\|u-w\|_\infty+c_2|R_p|$ for $c_2=(\pi(|a_0|+1)+4) e^{2|\Re a_0|}/(2 N^2)$.  In case $\Re a_0>0$, this is the best we can do.  Otherwise, by the previous paragraph we note $\|y-p\|_\infty \le 2 \|u-w\|_\infty + \min\{c_2,\pi/(2N)\} |R_p|$.

\bigskip

\centerline{\textsc{Appendix B: Compactness of the monodromy operator}}
\medskip

In theorem \ref{thm:main} we use theorem \ref{thm:apostivp} and an eigenvalue perturbation technique to estimate the difference between the eigenvalues of $U$ (eigenvalues of $\hU$, actually, but they are the same) and the eigenvalues of a matrix approximation to $\hU$, as the latter can be found numerically.  Of course, $\hU$ is not a matrix and its numerical approximation will (most naturally, anyway) act on a space of finite dimension.  We do not address, in theorem \ref{thm:main}, values in the spectrum of $\hU$ that are not eigenvalues.  However, as we show here, $U$ and $\hU$ are compact and so the only (possible) non-eigenvalue spectrum is $0\in\CC$.

\begin{defn}  Let $\Bcal$ be a Banach space with norm $\|\cdot\|$ and denote by $\Lop(\Bcal)$ the set of bounded operators on $\Bcal$.  We say $K\in\Lop(\Bcal)$ is \emph{compact} if the image $KS$ of any bounded subset $S\subset \Bcal$ is precompact, that is, if every sequence in $KS$ has a convergent subsequence.
\end{defn}

See section 21 of \cite{Lax} for a proof of the following proposition.

\begin{prop} The set of compact operators is a subspace of $\Lop(\Bcal)$ and forms a two-sided ideal (that is, if $L\in\Lop(\Bcal)$ and $K$ is compact then both $KL$ and $LK$ are compact).  Operators with finite rank are compact.  The spectrum $\sigma(K)$ of a compact operator $K$ consists of a countable set of complex numbers $\lam_k\in \CC$.  The only possible accumulation point of $\sigma(K)$ is zero; if $\dim \Bcal=\infty$ then $0\in \sigma(K)$. Nonzero elements of $\sigma(K)$ are eigenvalues.\end{prop}

\begin{thm*}  $U\in\Lop(C(I))$, defined by equation \eqref{Udef}, is compact.\end{thm*}

\begin{proof}  Section 22.2 of \cite{Lax} shows that integral operators $(N_k f)(t)=\int_{-1}^t k(t,s) f(s)\,ds$ are compact on $C(I)$ if the kernel $k(t,s)$ is continuous in each variable; this fact is an easy corollary of the Arzela-Ascoli theorem.  Note that multiplication of $f$ by the bounded functions $b,e^{\pm a(t+1)}$ are each bounded operators; denote by $M_b,M_\pm$.  The operator $\delta_1=f\mapsto f(1)$ is finite rank and thus compact.  Let $(N_1 f)(t)=\int_{-1}^t f(s)\,ds$, a compact integral operator.  It follows that $U=M_+ \left(\delta_1 + N_1 M_- M_b\right)$ is compact (because compact operators are a subspace and a two-sided ideal).
\end{proof}

\begin{cor*} $\hU\in\Lop({H^1})$ is compact. \end{cor*}

\begin{proof} Recall $\hU=U\circ i_{{H^1}}$ on ${H^1}$ and that the injection $i_{{H^1}}:{H^1} \hookrightarrow C(I)$ is bounded. \end{proof}

\bigskip

\centerline{\textsc{Appendix C: The Bauer-Fike theorem for matrices}}
\medskip

Let $\|\cdot\|$ be a norm for $n\times n$ matrices induced by a vector norm on $\CC^n$.  For such a norm we define $\cond(V) \equiv \|V\| \|V^{-1}\|$, the \emph{matrix condition number} of $V$.

\begin{thm*}[Bauer-Fike \cite{BauerFike}]  \label{thm:BF}  Suppose an $n\times n$ matrix $A$ is diagonalized by $V$, that is, suppose $A=V\Lambda V^{-1}$ for an invertible $n\times n$ matrix $V$ and $\Lambda=\diag(\lam_1,\dots,\lam_n)$, $\lam_i\in\CC$.  If $B$ is another $n\times n$ matrix with eigenvalue $\mu\in\CC$ then
\begin{equation}\label{BFest}
\min_{i=1,\dots,n} \left|\mu-\lam_i\right| \le \cond(V) \|B-A\|.
\end{equation}\end{thm*}

\begin{proof}  Apply theorem \ref{thm:BFH} supposing $Bx=\mu x$ for $x\in\CC^n$ satisfying $\|x\|=1$. \end{proof}

\begin{figure}[ht]
\widefigure{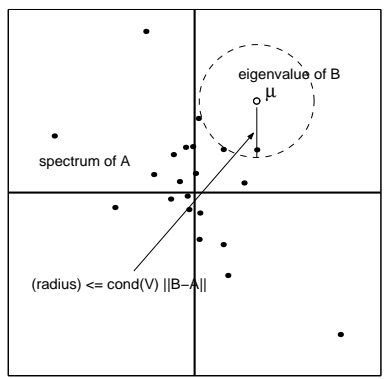}{2}{2}
\caption{The Bauer-Fike theorem bounds the distance between eigenvalues of $B$ and those of $A$ in terms of $\|B-A\|$ and the matrix condition number of a matrix of eigenvectors for $A$.} \label{fig:ABeigs}
\end{figure}
 
That is, if $A$ and $B$ are close in norm and $A$ is diagonalizable then any eigenvalue of $B$ will be close to one of the eigenvalues of $A$ if additionally $V$ is well-conditioned.   More technically, the Bauer-Fike theorem says that an upper bound for the condition number with respect to $\|\cdot\|$ of the eigenvalue problem for diagonalizable $A$ is given by the matrix condition number of any matrix of eigenvectors $V$.  Note that different diagonalizations of $A$ produce different values for $\cond(V)$, in general.  We suggest the reader explore this theorem in the following exercises.

\begin{exer}  Show that if $A$ is normal then $V$ can be chosen with $\cond(V)=1$.  (Thus the eigenvalue problem for normal matrices is well-conditioned.) \end{exer}
 
\begin{exer}  Let $A_1=\diag(1,2,3)$ and $A_2=[1\; 5 \; 10;\; 0\; 2 \; -\!\!10;\; 0\; 0\; 3]$ (\Matlab~notation).  These matrices have the same eigenvalues.  Show (numerically), by choosing random but small perturbations of the entries, that the eigenvalues of $A_2$ are much more sensitive to perturbations than those of $A_1$.\end{exer}

\begin{exer}  Let $A=A_2$ in the above exercise.  Let $\Delta A$ be a matrix with entries randomly chosen from an $N(0,1)$ distribution. Let $B=A+0.01 \Delta A$. Experiment numerically with different diagonalizations $V$ of $A$ and different matrix norms to make \eqref{BFest} as close to equality as possible. \end{exer}

Theorem \ref{thm:BFH} generalizes Bauer-Fike to Hilbert spaces, of course, but also includes the idea that $B,A$ do not have to be close in norm to get the result; they need only \emph{do close to the same thing} to $x$ so that $\|(B-A)x\|$ is small.

\end{document}